\pdfoutput=1 
\documentclass[11pt,reqno]{amsart}

\usepackage{amsmath, amssymb, graphicx, epsfig, verbatim}
\usepackage{epic}

\usepackage[margin=1in]{geometry}
\usepackage{booktabs, url}
\usepackage[section]{placeins}
\usepackage{tikz}
\usetikzlibrary{arrows.meta}

\newtheorem*{theoremA}{Theorem A}
\newtheorem*{theoremB}{Theorem B}

\theoremstyle{definition}

\theoremstyle{remark}

\newcommand{\PP}{\mathbf P}

\newcommand{\ZZ}{\mathbf Z}
\newcommand{\CC}{\mathbf C}

\newcommand{\Mod}{\operatorname{Mod}}
\newcommand{\panel}[3]{\begin{minipage}[t]{#1\textwidth}\centering\includegraphics[width=\linewidth]{#2}\par\smallskip{\small #3}\end{minipage}}

\numberwithin{equation}{section}
\graphicspath{{figures/}}

\begin{document}

\title{The monodromy of Xiao's genus-two fibrations in degrees 3, 4 and 5}
\author{Zolt\'an Szab\'o}
\address{Department of Mathematics, Princeton University, Princeton, NJ 08544, \
USA}
\email{szabo@math.princeton.edu}

\date{\today}

\begin{abstract}
  The goal of this paper is to determine the vanishing cycles of Xiao's genus-two fibrations $X_d\to\PP^1$, $d=3,4,5$. The fibers admit degree-$d$ maps to a
  fixed elliptic curve $E$. We subdivide the base $\PP^1$ by triangles and realize the surface $X_d$ as a degree-$d$ branched cover of $\PP^1\times E$.
This tiling is given by a topological classification of degree $d$ branched covers from a genus $2$ surface to a torus. 
The covering description gives an explicit algorithm for the vanishing cycles. A vanishing path is recorded by the sides it crosses in the tiling; each crossing changes a labelled picture of the fiber by a local \emph{hexagon move}, and the vanishing cycle is read from
the terminal picture. We obtain the $7$, $13$, and $31$ vanishing cycles as explicit curves in a marked reference fiber.
For $d=4,5$, these give new positive factorizations of the identity in the genus-two mapping class group, of types $(6,7)$ and $(12,19)$.
\end{abstract}
\maketitle

\section{Introduction}\label{sec:intro}

In~\cite{Xiao}  Xiao constructed genus-two fibrations whose fibers admit maps of degree $d$ to a fixed elliptic curve $E$.
In this paper we study the fibrations
\[
 S(E,d)\longrightarrow X(d),\qquad d=3,4,5,
\]
whose base $X(d)$ is isomorphic to $\PP^1$. Their types, the numbers of nonseparating and separating vanishing cycles, are respectively
\[
 (4,3),\qquad (6,7),\qquad (12,19),
\]
were also  computed by Xiao in~\cite{Xiao}. In this paper we  give explicit descriptions of these vanishing cycles and their monodromy factorizations.

The case $d=2$ was worked out by Matsumoto, whose relation $(t_{B_0}t_{B_1}t_{B_2}t_C)^2=1$ gives its monodromy of type $(6,2)$~\cite{Mat}.
For $d=3$, Baykur and Korkmaz reverse-engineered a positive factorization of type $(4,3)$~\cite{BK}.
Its Hurwitz equivalence with the monodromy of Xiao's fibration remained open~\cite[Remark~5]{Nak}. Such factorizations provide building blocks for constructing exotic smooth four-manifolds and further genus-two fibrations; see, for example~\cite{AM,BK,Nak,OS} for some really interesting applications.

There are more than one way to study these fibrations. For example $d=3$ case can also be directly studied using Xiao's plane model, an alternative
description of $S(E,3)$ using a double branched cover, see \cite{Xiao} and also \cite{AM}.
The $d=3$ vanishing cycles were computed by the author using Xiao's plane model.
This work was presented at the ICTP summer school in Trieste in June 2025~\cite{Trieste}. Huang~\cite{Huang} used this computation to prove the Hurwitz equivalence with the Baykur--Korkmaz factorization. An account of the degree-three calculation, together with the elliptic-covering approach discussed in this paper, was presented at the Simons Collaboration annual meeting in March 2026~\cite{Simons}.
For a later account of the $d=3$ case using Xiao's plane model,
see~\cite{Akh}.

The goal of this paper is to give a detailed description of the elliptic-covering approach and results presented in~\cite{Simons}. The method itself works the same way for $d=3,4,5$. In this paper we give the detailed calculation in degrees three and four, and also give the explicit monodromy factorization in degree five.

\subsection{Xiao's fibrations}\label{subsec:xiao-universal}
Fix $d>2$ and write $E=\CC/(\ZZ+\ZZ\tau_0)$. For each $\tau$ in the upper half-plane, Xiao writes down a lattice in $\CC^2$ and a genus-two curve $C_\tau$ whose Jacobian is the quotient of $\CC^2$ by this lattice. The curve comes with a degree-$d$ map $C_\tau\to E$ that is primitive, meaning that it does not factor through a nontrivial isogeny. Note that changing $\tau$ by a matrix in $\mathrm{SL}_2(\ZZ)$ that is congruent to the identity modulo~$d$ gives an isomorphic cover. The quotient of the upper half-plane by the group of these matrices, compactified by adding its cusps, is the curve $X(d)$. For $d=3,4,5$ it is well known that $X(d)$ is a sphere and that its $4$, $6$, $12$ cusps correspond to the vertices of a regular tetrahedron, octahedron, and icosahedron.

Xiao's family $S(E,d)\to X(d)$ is universal for normalized primitive covers $C\to E$~\cite[Theorem~3.10, Corollary~2]{Xiao}; see also~\cite{Kaya,KaniHurwitz}. For $d=3,4,5$, $S(E,d)$ has $7$, $13$, and $31$ singular fibers, each of which is either an irreducible curve with one node or two elliptic curves meeting transversally in one point~\cite[Lemma~3.11, Theorem~3.16, p.~52]{Xiao}, so it is a genus-two Lefschetz fibration.

We identify the tiled base with $X(d)$ by moving the branch values of a reference cover. Independently, we construct the corresponding branched cover of $\PP^1\times E$ and identify the resulting fibration using Xiao's universality and the singular-fiber count.

We write $f\colon X\to B$ for the family constructed below. For the explicit calculations we take
\begin{equation}\label{eq:E}
 E=\CC/(\ZZ+\ZZ\omega),\qquad \omega=e^{2\pi i/3},
\end{equation}
whose order-three symmetry simplifies the pictures. The extension to arbitrary elliptic curves is discussed in Section~\ref{subsec:other-E}.

\subsection{The topological construction and the algorithm}\label{subsec:intro-construction}
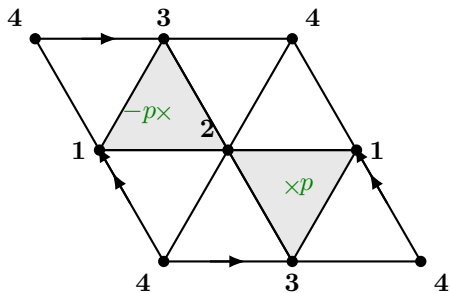
\begin{figure}[!htbp]
 \centering
 \begin{tikzpicture}[scale=3.4,every node/.style={font=\small}]
  
  \coordinate (A) at (0,0); \coordinate (B) at (1,0); \coordinate (C) at (0.5,0.8660); \coordinate (D) at (-0.5,0.8660);
  \coordinate (MAB) at (0.5,0); \coordinate (MBC) at (0.75,0.4330); \coordinate (MCD) at (0,0.8660); \coordinate (MDA) at (-0.25,0.4330); \coordinate (O) at (0.25,0.4330);
  
  \fill[gray!18] (MAB)--(MBC)--(O)--cycle; \fill[gray!18] (MDA)--(MCD)--(O)--cycle;
  
  \draw[thick] (A)--(B)--(C)--(D)--cycle;
  \draw[thick] (A)--(C);
  \draw[thick] (MAB)--(MBC)--(O)--cycle; \draw[thick] (MDA)--(MCD)--(O)--cycle;
  \draw[thick] (MAB)--(O)--(MCD); \draw[thick] (MBC)--(O)--(MDA);
  
  \draw[thick,-{Latex[length=2.2mm]}] (0.18,0)--(0.32,0); \draw[thick,-{Latex[length=2.2mm]}] (-0.32,0.8660)--(-0.18,0.8660);
  \draw[thick,-{Latex[length=2.2mm]}] (-0.13,0.2252)--(-0.20,0.3464); \draw[thick,-{Latex[length=2.2mm]}] (-0.19,0.3291)--(-0.26,0.4503);
  \draw[thick,-{Latex[length=2.2mm]}] (0.87,0.2252)--(0.80,0.3464); \draw[thick,-{Latex[length=2.2mm]}] (0.81,0.3291)--(0.74,0.4503);
  
  \node[green!50!black] at (0.50,0.289) {$\times$}; \node[green!50!black,right=-1pt] at (0.50,0.289) {$p$};
  \node[green!50!black] at (0.00,0.577) {$\times$}; \node[green!50!black,left=-1pt] at (0.00,0.577) {$-p$};
  
  \foreach \P/\L/\pos in {A/4/below left, B/4/below right, C/4/above right, D/4/above left, MAB/3/below, MCD/3/above, MBC/1/right, MDA/1/left, O/2/above left}
   { \fill (\P) circle (0.022); \node[\pos=1pt,font=\small\bfseries] at (\P) {\L}; }
 \end{tikzpicture}
 \caption{The elliptic curve $E=\CC/(\ZZ+\ZZ\omega)$ as a parallelogram with opposite sides identified, subdivided into eight triangles by the four points of $E[2]$: $4=0$, $3=\tfrac12$, $1=\tfrac\omega2$, and $2=\tfrac{1+\omega}2$. The involution $-1$ interchanges the two triangles over each face of the tetrahedron on $T$. The shaded triangles lie over the face $123$ and contain the marked branch values $\pm p$.}
 \label{fig:intro-torus}
\end{figure}

\begin{figure}[!htbp]
 \centering
 \panel{.27}{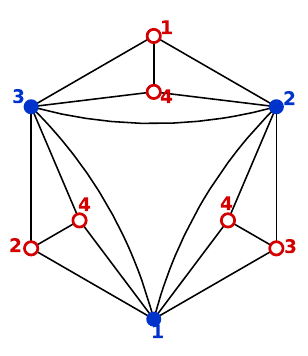}{(a) $d=3$.}\hfill
 \panel{.27}{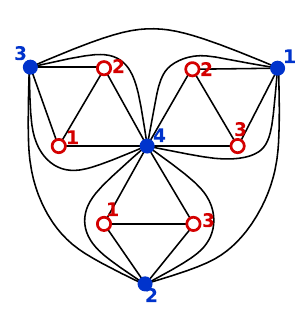}{(b) $d=4$.}\hfill
 \panel{.27}{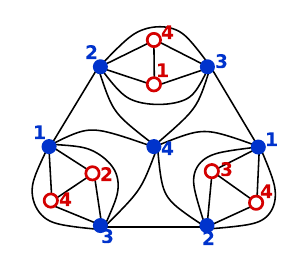}{(c) $d=5$.}
 \caption{The reference coverings for $d=3,4,5$, drawn on the six-point sphere $H$. Each vertex is labelled $1$, $2$, $3$ or $4$ by its image in $T$, with $4=q(0)$; the red vertices are the images of the Weierstrass points. In (a) and (c), the exceptional hexagon is the outer face. In (b), the two branch values coincide at $0$, so there is no hexagon and the central blue vertex has local degree four; the two ramification points remain distinct in $C$.}
 \label{fig:intro-central}
\end{figure}

Let $h\colon C\to E$ be a degree-$d$ cover by a genus-two curve, and let $\iota$ be the hyperelliptic involution of $C$. After a translation of $E$ we may assume $h\circ\iota=-h$; then the six Weierstrass points of $C$ map to the two-torsion points $E[2]$, and the two simple branch values of $h$ are $\pm p$. Up to isomorphism, the cover is determined by its monodromy $\pi_1(E\setminus\{p,-p\})\to S_d$, and we encode this monodromy in a picture. The map $h$ descends to a degree-$d$ map from the six-point sphere $H=C/\langle\iota\rangle$ to the four-point sphere $T=E/\{\pm1\}$, whose marked points are the images $q(E[2])$ under the quotient map $q\colon E\to T$, labelled $1,2,3,4$ as in Figure~\ref{fig:intro-torus}. Joining these points by six edges makes $T$ a tetrahedron, and pulling back its faces gives a picture of $H$ consisting of triangles and one exceptional hexagon, the preimage of the face containing $q(p)$, with the six Weierstrass images marked red (Figure~\ref{fig:intro-central}). The labelled picture determines the monodromy, so two covers are isomorphic exactly when their pictures are.

\begin{figure}[p]
 \centering
 \panel{.35}{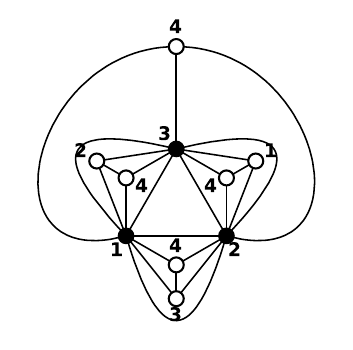}{(a) $d=3$.}\hfill
 \panel{.52}{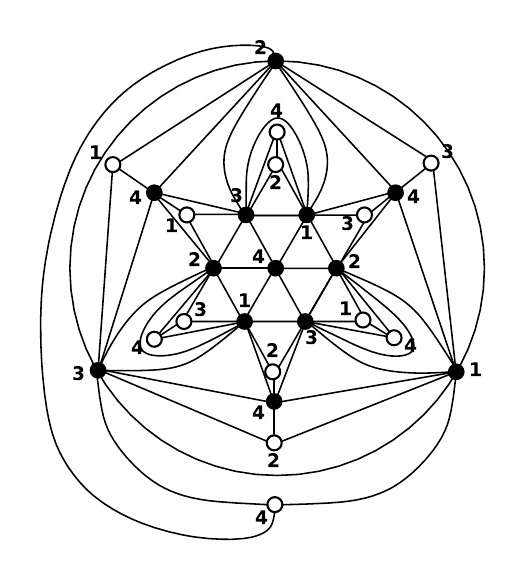}{(b) $d=4$.}
 \par\smallskip
 \panel{.44}{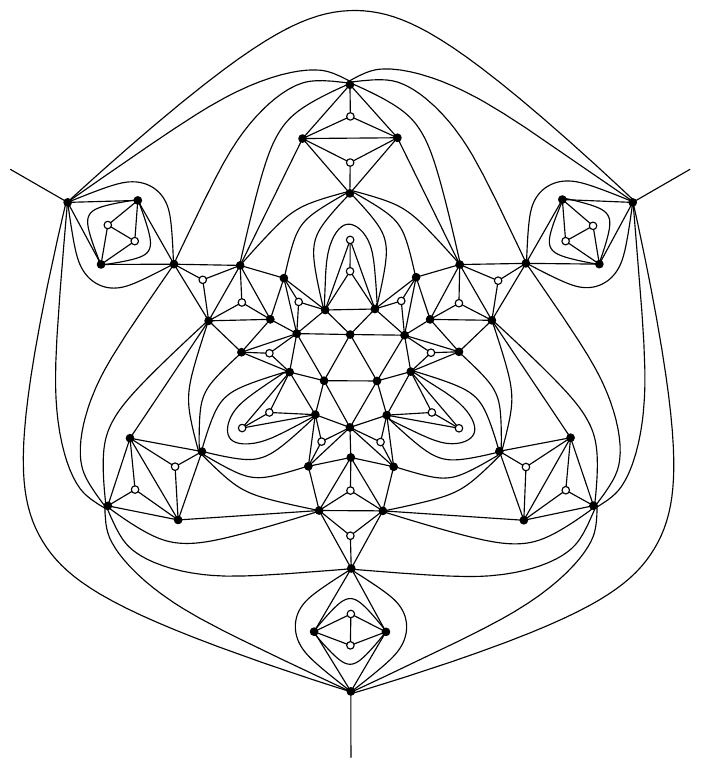}{(c) $d=5$.}
 \caption{The base $B$ of the fibration, tiled by triangles whose vertices are labelled by the four branch points of $q$. The hollow circles are the vertices with three edges, which are the critical values; the black dots are the other vertices. The base point $b_0$ is the center of the central triangle in (a) and (c) and the central vertex in (b). In (a) the singular vertex $P$ at the top is the point at infinity; in (c) the three lines leaving the picture meet at the singular vertex at infinity.}
 \label{fig:intro-tilings}
\end{figure}

Moving the branch values $\pm p$ gives the map $G\colon B\to T$, with $G(t)=q(p_t)$. Pulling back the tetrahedron by $G$ tiles $B$ by $16,48,160$ triangles (Figure~\ref{fig:intro-tilings}). Crossing a side changes the marked covering picture by a hexagon move. At a nodal endpoint, the vanishing cycle is the lift of an arc or circle in the terminal picture. Sections~\ref{sec:move}--\ref{subsec:d5} describe the move and carry out the three calculations.

\subsection{Results}\label{subsec:results}
Give $B\cong\PP^1$ the complex structure for which $G\colon B\to T$ is holomorphic. The corresponding degree-$d$ cover
\[
 \Pi=(f,g)\colon X\longrightarrow B\times E
\]
is branched along
\[
 \Delta=\{(t,e):G(t)=q(e)\}.
\]

\begin{samepage}
\begin{theoremA}
For $d=3,4,5$, Xiao's surface $S(E,d)$ for the hexagonal curve $E$ is the $d$-fold branched covering $X$ of $\PP^1\times E$ just described, with the complex structure for which $\Pi$ is holomorphic, and Xiao's fibration $S(E,d)\to X(d)$ is $f$; in particular $X(d)=B\cong\PP^1$. The branch curve $\Delta$ has $0$, $4$, $24$ nodes and $3$, $9$, $27$ cusps, over the vertices with six and nine edges, and the $7$, $13$, $31$ singular fibers of $f$ lie over its simple tangencies with the vertical fibers, at the hollow circles.
\end{theoremA}
\end{samepage}

The proof, in Section~\ref{sec:surface}, uses the fact that every fiber of $f$ carries a degree-$d$ map to the fixed curve $E$, so $f$ is a fibration of type $(E,d)$ in Xiao's sense, and Xiao's uniqueness theorem for such fibrations identifies it with $S(E,d)$ once the number of singular fibers is known. The remaining results are theorems about $f$ and are proved topologically.

\begin{theoremB}
Let $d\in\{3,4,5\}$, let $b_0$ be the central point of $B=X(d)$, and let $\gamma_1,\dots,\gamma_n$, $n=7,13,31$, be the system of vanishing paths of Figures~\ref{fig:d3-paths}, \ref{fig:paths} and~\ref{fig:d5-paths}, numbered as in~\eqref{eq:d3order}, \eqref{eq:d4order} and~\eqref{eq:d5order}. The vanishing cycles $\delta_1,\dots,\delta_n$ of Xiao's fibration $S(E,d)\to X(d)$ along these paths are the curves in the fiber $\Sigma=f^{-1}(b_0)$ shown in Figures~\ref{fig:d3-summary}, \ref{fig:all-cycles} and~\ref{fig:d5-summary}; of these, $4$, $6$, $12$ are nonseparating and $3$, $7$, $19$ are separating. The monodromy factorization of Xiao's fibration is
\[
 t_{\delta_n}\circ\cdots\circ t_{\delta_1}=1
\]
in $\Mod(\Sigma_2)$, where $t_\delta$ denotes the right-handed Dehn twist. On the six-point sphere, this gives the corresponding identity in half-twists and squared twists.
\end{theoremB}

The surfaces $X$ have $(e,\sigma)=(3,-3),(9,-5),(27,-11)$, respectively~\cite{Mat,En}. Algebraic intersections of vanishing cycles are divisible by $d$.

The same method would also work, in principle, for $d>5$. For example, for $d=6$ the base is a torus, and it would be interesting to compute its tiling. However, as the calculations in this paper show, the vanishing cycles get progressively more complicated as the singular fiber gets farther from the reference fiber in the tiling, that is, as the vanishing path crosses more sides.

\subsection*{Acknowledgments}
The author would like to thank \.{I}nan\c{c} Baykur for suggesting the study of Xiao's genus-two fibrations with $d>3$, for sharing his computation of the fundamental group in degree four, and for many helpful discussions throughout this project.

The author is also grateful to Eaman Eftekhary, Evan Huang, J\'anos Koll\'ar, Peter Kronheimer, Ciprian Manolescu, Tomasz Mrowka, Burak Ozbagci, Peter Ozsv\'ath, and Andr\'as Stipsicz for many helpful conversations. The author was partially supported by the Simons Collaboration Grant on New Structures in Low-Dimensional Topology.

{\it AI Statement:}
In the course of the preparation of this paper OpenAI's ChatGPT and Anthropic's Claude were used for typesetting, language editing,
the creation of some of the figures and writing code to do some additional checks.
The author assumes full intellectual responsibility for the mathematical statements, computations, and any remaining errors.

\section{Coverings of the torus and their pictures}\label{sec:covers}

This section fixes notation for the covers and pictures of Section~\ref{subsec:intro-construction}. Throughout, $E$ and its labeling
are as in Figure~\ref{fig:intro-torus}. Let $h\colon C\to E$ be a degree-$d$ branched cover by a genus-two surface $C$, normalized as in Section~\ref{subsec:intro-construction} so that
\begin{equation}\label{eq:equivariance}
 h\circ\iota=-h.
\end{equation}
The induced map $\bar h\colon H=C/\langle\iota\rangle\to T$ satisfies $q\circ h=\bar h\circ\pi$, where $\pi\colon C\to H$ is the quotient map. We call the cover generic when its two simple branch values $p,-p$ are distinct; the remaining cases are treated in Section~\ref{subsec:corners}. For the pictures we also assume that $q(p)$ lies in the interior of a face; crossing a side is the hexagon move of Section~\ref{subsec:crossing}.

Pulling back the tetrahedron gives the picture of $H$: its red vertices are the Weierstrass images, and its blue vertices are the remaining preimages of $q(E[2])$. The map $\bar h$ has local degree one at red vertices and two at blue vertices. We normalize further, by a translation by a point of $E[2]$, so that three Weierstrass points map to $0$ when $d$ is odd and none do when $d=4$. In the reference pictures the Weierstrass distribution, that is, the numbers of red vertices over $1,2,3,4$, is then $(1,1,1,3)$ for odd $d$ and $(2,2,2,0)$ for $d=4$; since $B$ is connected, the same holds for every cover in the family. The exceptional hexagon lies over the face containing $q(p)$ and has cyclic labels $i,j,k,i,j,k$.

The covering is determined by its monodromy
\begin{equation}\label{eq:monodromy}
 \mu\colon\pi_1\bigl(E\setminus\{p,-p\},e_*\bigr)\longrightarrow S_d,
\end{equation}
up to conjugacy. Two covers are isomorphic exactly when their labelled pictures are isomorphic. Here an isomorphism of pictures is an orientation-preserving homeomorphism preserving the cells, labels and colors. The $d$ preimages in $H$ of a point of $T$ correspond to the $d$ sheets of $h$ over either of its preimages in $E$, so $\mu$ can be read off by lifting loops through the triangles of the picture. The reference pictures are Figures~\ref{fig:intro-central}(a) and~(c) for $d=3,5$. For $d=4$ the reference cover has no exceptional hexagon, we use the perturbed picture of Figure~\ref{fig:first-move}(a) (Section~\ref{subsec:d4-reference}).

Write $S$ for the sphere $H$ of the reference cover, $W\subset S$ for its six red points, and $\pi\colon\Sigma\to S$ for the double cover branched along $W$; thus $\Sigma=f^{-1}(b_0)$ is the reference fiber. In degrees three and five, $r_i$ is the red point over $i=1,2,3$. The red points over $4$ adjacent to $r_1,r_2,r_3$ are $a,b,c$ in degree three and $b,c,a$ in degree five. In degree four, the red points of Figure~\ref{fig:first-move}(a), counterclockwise from the leftmost, are $r_1,r'_1,r_3,r'_3,r_2,r'_2$. The order-three rotation of the reference picture (relabelling $1\mapsto2\mapsto3\mapsto1$; see Section~\ref{subsec:symmetry}) permutes each of the triples $r_1,r_2,r_3$; $a,b,c$; and $r'_1,r'_2,r'_3$ cyclically.

\subsection*{Conventions}
The drawings have the standard orientation, with $1\to2\to3$ positive around face $123$. In a product, the rightmost factor acts first. Write $t_\delta$ for the right-handed Dehn twist about a simple closed curve $\delta\subset\Sigma$, acting on homology by $t_\delta(x)=x+(\delta\cdot x)\delta$. Write $H_\alpha$ for the right-handed half-twist along an arc $\alpha\subset S$ joining two points of $W$, and $T_\gamma$ for the right-handed Dehn twist about a simple closed curve $\gamma\subset S\setminus W$.

\section{The hexagon move and the tilings of the base}\label{sec:move}

\subsection{The hexagon move}\label{subsec:crossing}\label{subsec:transport}
\begin{figure}[htbp]
 \centering
 \begin{tikzpicture}[scale=1.15,every node/.style={font=\small}]
 \begin{scope}
  \coordinate (i1) at (1.3,0); \coordinate (k1) at (0.92,0.92); \coordinate (j1) at (0,1.3); \coordinate (l1) at (-0.92,0.92);
  \coordinate (i2) at (-1.3,0); \coordinate (k2) at (-0.92,-0.92); \coordinate (j2) at (0,-1.3); \coordinate (l2) at (0.92,-0.92);
  \fill[gray!25] (i1)--(k1)--(j1)--(i2)--(k2)--(j2)--cycle;
  \draw[thick] (i1)--(k1)--(j1)--(l1)--(i2)--(k2)--(j2)--(l2)--cycle;
  \draw[thick] (j1)--(i2); \draw[thick] (j2)--(i1);
  \node[above right] at (i1) {$i_1$}; \node[above right] at (k1) {$k_1$}; \node[above] at (j1) {$j_1$}; \node[above left] at (l1) {$l_1$};
  \node[left] at (i2) {$i_2$}; \node[below left] at (k2) {$k_2$}; \node[below] at (j2) {$j_2$}; \node[below right] at (l2) {$l_2$};
  \node at (0,0) {$\bullet$}; \node[right] at (0.05,0.02) {$P$};
  \node at (0,-2) {before: $p\in\Delta_1$};
 \end{scope}
 \begin{scope}[xshift=4.6cm]
  \coordinate (i1) at (1.3,0); \coordinate (k1) at (0.92,0.92); \coordinate (j1) at (0,1.3); \coordinate (l1) at (-0.92,0.92);
  \coordinate (i2) at (-1.3,0); \coordinate (k2) at (-0.92,-0.92); \coordinate (j2) at (0,-1.3); \coordinate (l2) at (0.92,-0.92);
  \fill[gray!25] (i1)--(j1)--(l1)--(i2)--(j2)--(l2)--cycle;
  \draw[thick] (i1)--(k1)--(j1)--(l1)--(i2)--(k2)--(j2)--(l2)--cycle;
  \draw[thick] (i1)--(j1); \draw[thick] (i2)--(j2);
  \node[above right] at (i1) {$i_1$}; \node[above right] at (k1) {$k_1$}; \node[above] at (j1) {$j_1$}; \node[above left] at (l1) {$l_1$};
  \node[left] at (i2) {$i_2$}; \node[below left] at (k2) {$k_2$}; \node[below] at (j2) {$j_2$}; \node[below right] at (l2) {$l_2$};
  \node at (0,0) {$\bullet$}; \node[right] at (0.05,0.02) {$P$};
  \node at (0,-2) {after: $p\in\Delta_2$};
 \end{scope}
 \end{tikzpicture}
 \caption{The hexagon move over $\Delta_1\cup\Delta_2$. The chords lift the crossed side $ij$; the shaded region is the exceptional hexagon. Here $h_1,\ldots,h_6=i_1,k_1,j_1,i_2,k_2,j_2$ and $x,y=l_2,l_1$.}
 \label{fig:disk}
\end{figure}
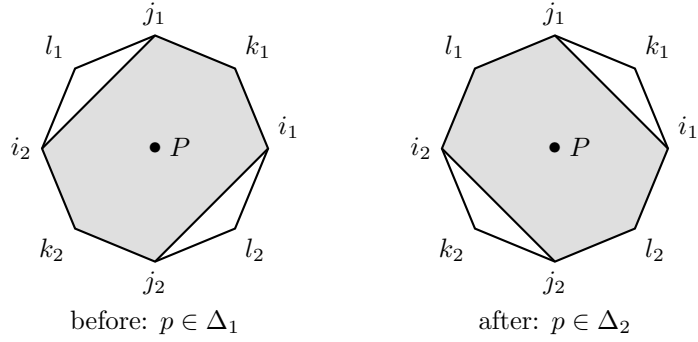

Suppose $p$ crosses $s=ij$ from $\Delta_1=ijk$ to $\Delta_2=ijl$. Label the corners of the old hexagon $h_1,\ldots,h_6$ so that the sides $h_6h_1$ and $h_3h_4$ lie over $s$, and let $x$ and $y$ be the vertices opposite these sides in the adjacent triangles. The crossing replaces $h_6h_1,h_3h_4$ by $h_1h_3,h_4h_6$, drawn inside the old hexagon alongside $h_1h_2h_3$ and $h_4h_5h_6$. The new hexagon has cyclic corners $x,h_1,h_3,y,h_4,h_6$ (Figure~\ref{fig:disk}).

A crossing word lists the crossed sides, such as $34,13,12$. Applying the corresponding moves gives the terminal picture in the reference sphere. When we transport a picture along a path, we keep its marking; when we enumerate covering classes, we identify pictures up to orientation-preserving labelled isomorphism.

\subsection{The tilings of the base}\label{subsec:base}
For each face of the tetrahedron, take one triangle for each covering class with $q(p)$ at a fixed interior point of that face, and glue the triangles along their sides according to the hexagon moves. Starting from the reference picture and filling in the vertices gives
\begin{equation}\label{eq:G}
 G\colon B\longrightarrow T.
\end{equation}
If a covering class returns to itself after $r$ turns around a vertex of $T$, the corresponding vertex of $B$ has local degree $r$ and valence $3r$. In degree four, start with the generic picture of Figure~\ref{fig:first-move}(a).

Enumerating the covers with the Weierstrass distributions of Section~\ref{sec:covers} and applying the three moves to each picture gives the spheres of Figure~\ref{fig:intro-tilings}. The degrees and ramification profiles of the maps $G$ are
\[
\begin{array}{c|c|cc}
 d&\deg G&\text{over each of }1,2,3&\text{over }4\\ \hline
 3&4&1^1\,3^1&1^4\\
 4&12&1^3\,3^3&1^4\,2^4\\
 5&40&1^6\,2^8\,3^6&1^{13}\,3^9
\end{array}
\]
where $r^a$ denotes $a$ points of local degree $r$. Thus the generic part of $B$ is Xiao's generic-cover locus; the local models of Section~\ref{sec:fibration} extend the identification across the vertices.

\subsection{The order-three symmetry}\label{subsec:symmetry}
Relabelling $1\mapsto2\mapsto3\mapsto1$ commutes with the hexagon move and induces the order-three rotation $\sigma$ of Figure~\ref{fig:intro-tilings}, with
\begin{equation}\label{eq:sigma}
 G\circ\sigma=\bar\omega\circ G,
\end{equation}
where $\bar\omega$ is induced by multiplication by $\omega$ on $E$. The corresponding rotation of the reference picture acts as in Section~\ref{sec:covers}. The fixed points of $\sigma$ are $b_0$ and the nodal vertex over $4$, denoted $P,M_0,\infty$ in degrees $3,4,5$.

The other critical values form orbits of size three. Rotated paths to these orbits, together with one path to the fixed critical value, reduce the calculation to $3,5,11$ transports.

\section{The fibration and its vanishing cycles}\label{sec:fibration}

\subsection{Collisions}\label{subsec:corners}
At a vertex of $B$ over $i$, the two branch values $\pm p$ come together at the point $e_i\in E[2]$ labelled $i$. Choose local coordinates $u$ on $B$ and $v$ on $E$, centered at the vertex and at $e_i$, such that $q(v)=v^2$ and $G(u)=u^r$, where $r$ is the local degree of $G$. The branch curve $\Delta$ of~\eqref{eq:Delta} is locally
\begin{equation}\label{eq:local-branch}
 \Delta\colon v^2=u^r.
\end{equation}
Just before the collision, the exceptional hexagon has two corners over $i$; they determine the local covering. The possibilities are the following. 
\begin{center}
\begin{tabular}{@{}clcll@{}}
\toprule
Case & $i$-corners & $r$ & fiber  \\
\midrule
(A) & distinct red vertices & 1 & irreducible, one node \\
(B) & the same blue vertex & 1 & two elliptic curves \\
(C) & one red, one blue & 3 &  smooth \\
(D) & distinct blue vertices & 2 & smooth \\
\bottomrule
\end{tabular}
\end{center}
The branch transpositions are equal in (A) and (B), share one letter in (C), and are disjoint in (D). In Figure~\ref{fig:intro-tilings}, the vertices with $r=1,2,3$ are those with three, six and nine edges.

\subsection{The total space}\label{subsec:total-space}
The covers $h_t\colon C_t\to E$, for $t\in B$ away from the vertices, form a family (Section~\ref{subsec:base}); completing it with these local models gives a degree-$d$ branched cover
\[
 \Pi=(f,g)\colon X\longrightarrow B\times E
\]
with branch curve
\begin{equation}\label{eq:Delta}
 \Delta=\{(t,e)\in B\times E:q(e)=G(t)\}.
\end{equation}
Since the local models are smooth, the surface $X$ is smooth; moreover $g|_{C_t}=h_t$, and $f\colon X\to B$ is a genus-two Lefschetz fibration with $7,13,31$ one-node fibers, respectively, at the vertices of local degree one (the hollow circles of Figure~\ref{fig:intro-tilings}). Write
\begin{equation}\label{eq:global-monodromy}
 \rho\colon\pi_1\bigl((B\times E)\setminus\Delta\bigr)\longrightarrow S_d
\end{equation}
for the covering monodromy; meridians of $\Delta$ act by transpositions.

\subsection{Reading the vanishing cycle}\label{subsec:terminal}
For a vanishing path ending at a critical value over $i$, join the two $i$-corners across the terminal hexagon. In case~(A) this gives an arc between red points, whose lift to $\Sigma$ is the nonseparating vanishing cycle. In case~(B) it gives a circle through the repeated blue vertex, separating the red points into two triples; its lift is the separating vanishing cycle. Case~(B) occurs at every critical value over $1,2,3$ for odd $d$, and over $4$ for $d=4$, since then there is only one red vertex over each of $1,2,3$, respectively none over $4$.

\section{The vanishing cycles in degree three}\label{sec:d3}\label{subsec:d3}

The seven vanishing cycles come from three transports, $R_1,Q_1,P$, and the rotations of the first two.

\subsection{The reference cover and paths}\label{subsec:d3-reference}
\begin{figure}[htbp]
 \centering
 \includegraphics[width=.48\textwidth]{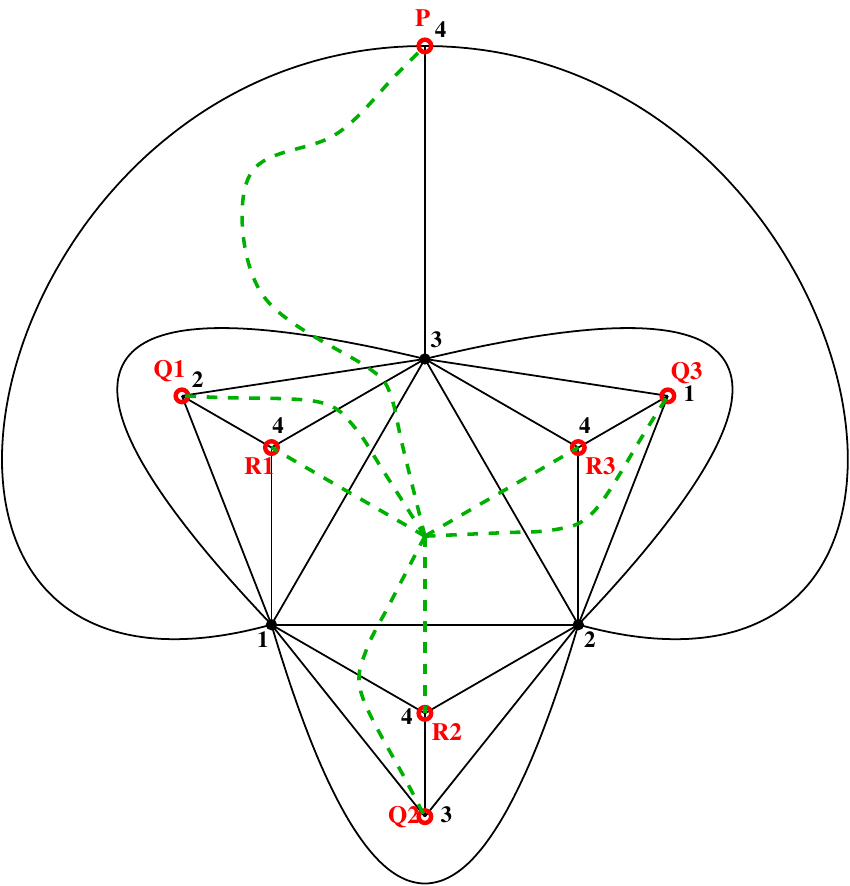}
 \caption{The degree-three vanishing paths. The base point $b_0$ lies in the central triangle; the red circles are the critical values.}
 \label{fig:d3-paths}
\end{figure}

Use the reference picture of Figure~\ref{fig:intro-central}(a) and the seven paths of Figure~\ref{fig:d3-paths}. Their crossing words are listed in Table~\ref{tab:d3-words}, and we take them in the order
\begin{equation}\label{eq:d3order}
 P,\ Q_1,\ R_1,\ Q_2,\ R_2,\ Q_3,\ R_3.
\end{equation}

\begin{table}[!htbp]
\centering\small
\begin{tabular}{@{}cccc@{}}
\toprule
path & crossing word & terminal label & terminal curve\\
\midrule
$R_1,\ R_2,\ R_3$ & $13;\quad 12;\quad 23$ & $4$ & arc\\
$Q_1,\ Q_2,\ Q_3$ & $13,34;\quad 12,14;\quad 23,24$ & $2;\ 3;\ 1$ & circle\\
$P$ & $13,34,23,13$ & $4$ & arc\\
\bottomrule
\end{tabular}
\medskip
\caption{The degree-three crossing words.}
\label{tab:d3-words}
\end{table}

\subsection{The three transports}\label{subsec:d3-transports}
\begin{figure}[htbp]
 \centering
 \panel{.40}{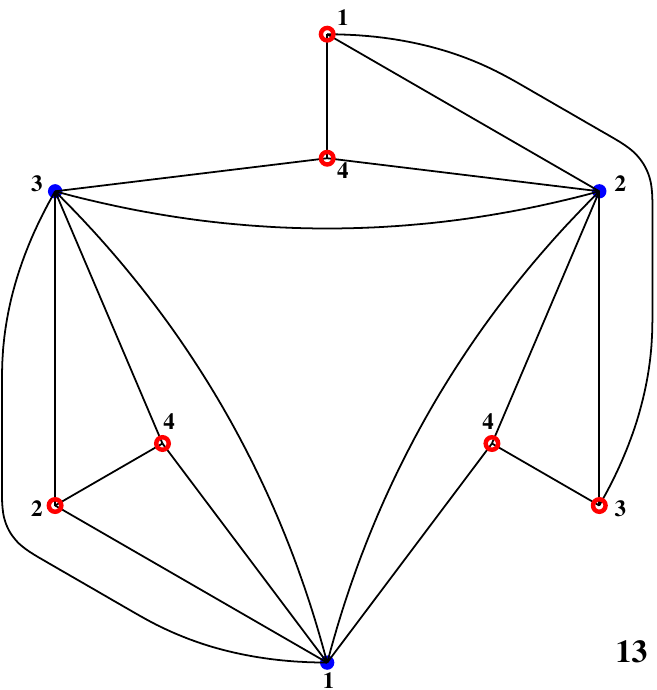}{(a) After $13$.}\hspace{.08\textwidth}
 \panel{.40}{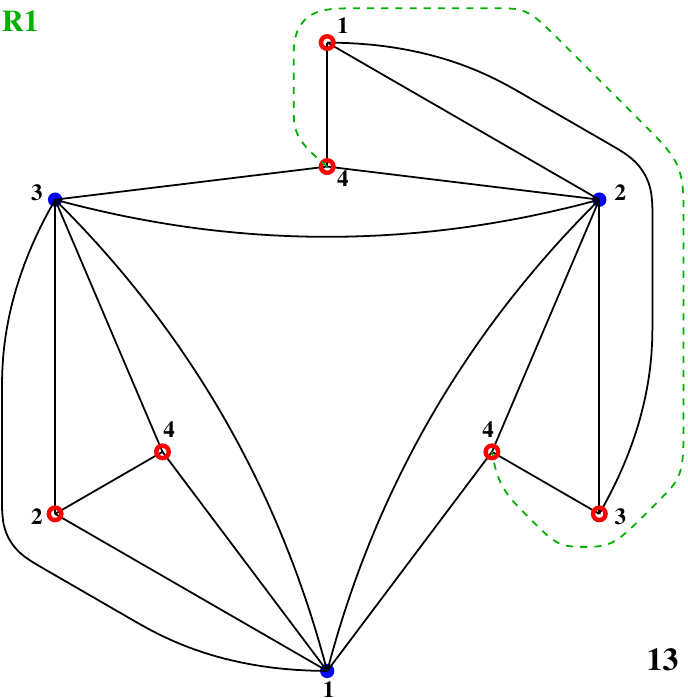}{(b) $R_1$: the arc from $a$ to $c$.}
 \par\medskip
 \panel{.40}{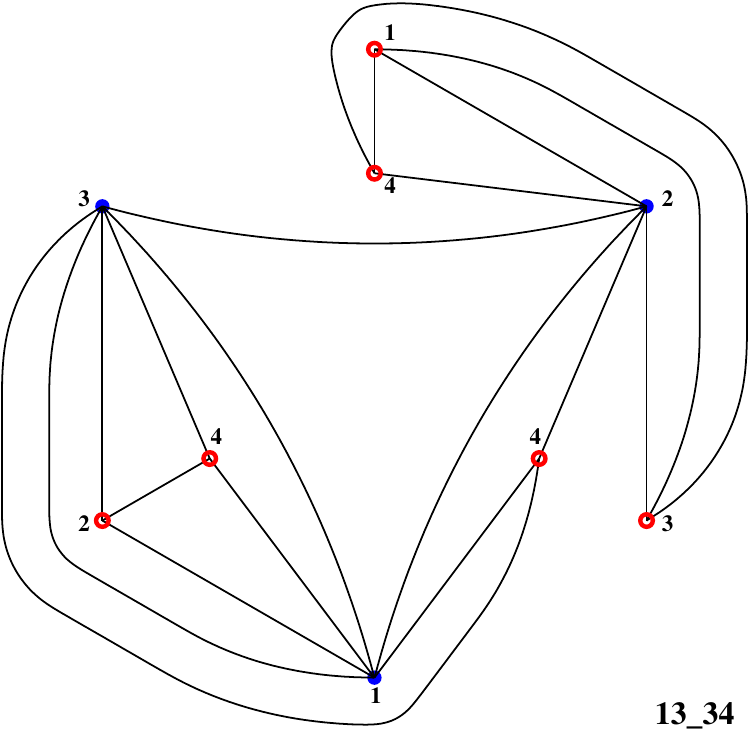}{(c) After $13,34$.}\hspace{.08\textwidth}
 \panel{.40}{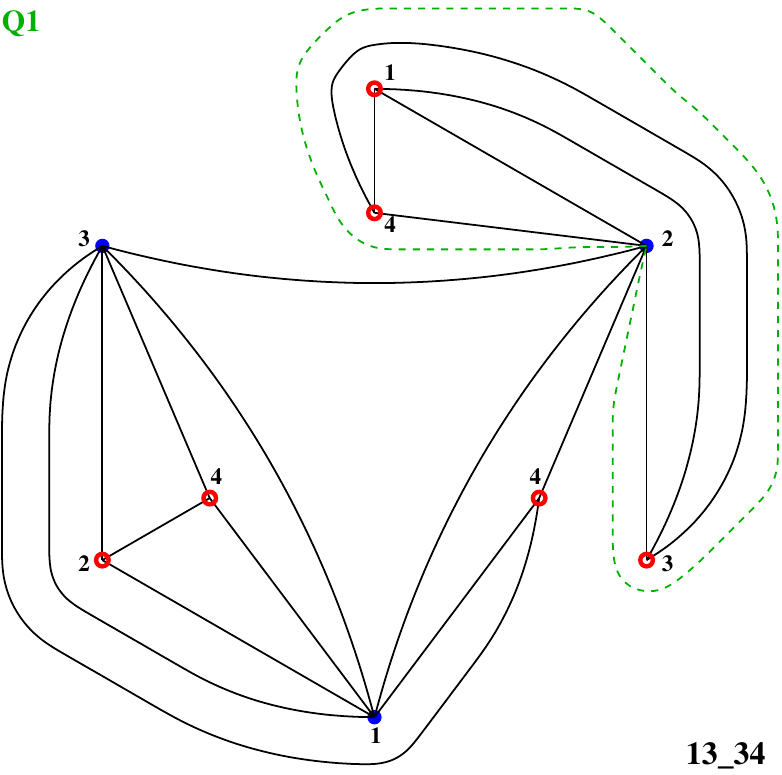}{(d) $Q_1$: the circle through the blue $2$.}
 \caption{The transports $R_1$ and $Q_1$: pictures after $13$ and $13,34$, with their terminal curves in green.}
 \label{fig:d3-r1q1}
\end{figure}

Figure~\ref{fig:d3-r1q1} shows the pictures after the moves $13$ and $13,34$. After $13$, the two $4$-corners are $a,c$, joined by the terminal arc $\alpha_{R_1}$. After $13,34$, the two $2$-corners are the same blue vertex; the terminal circle $\gamma_{Q_1}$ separates $\{r_1,r_3,a\}$ from $\{r_2,b,c\}$.

\begin{figure}[htbp]
 \centering
 \panel{.47}{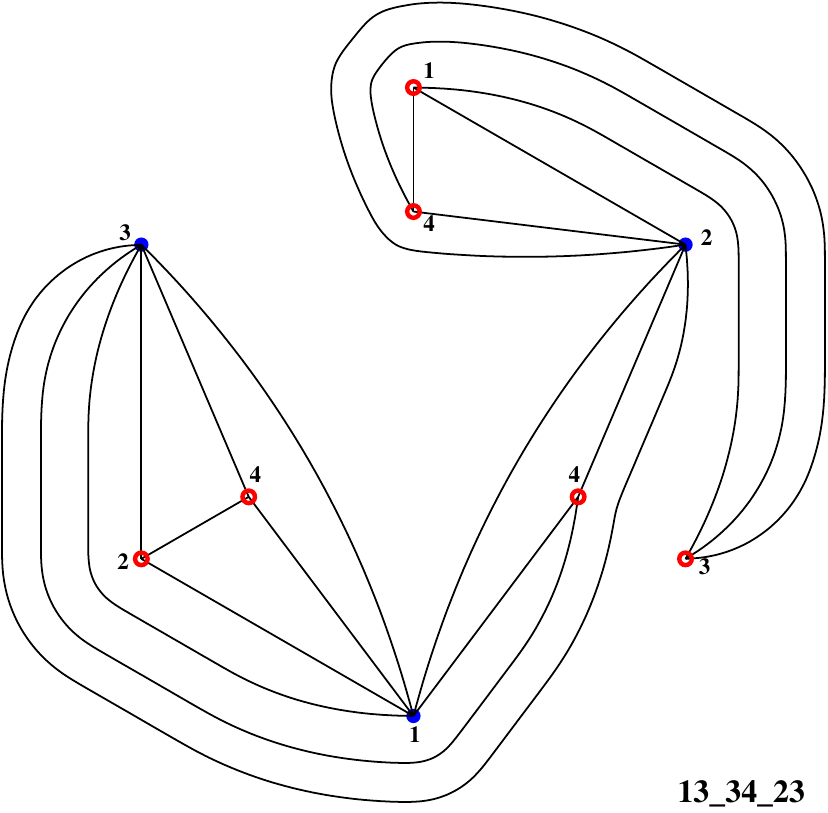}{(a) After $13,34,23$.}\hfill
 \panel{.47}{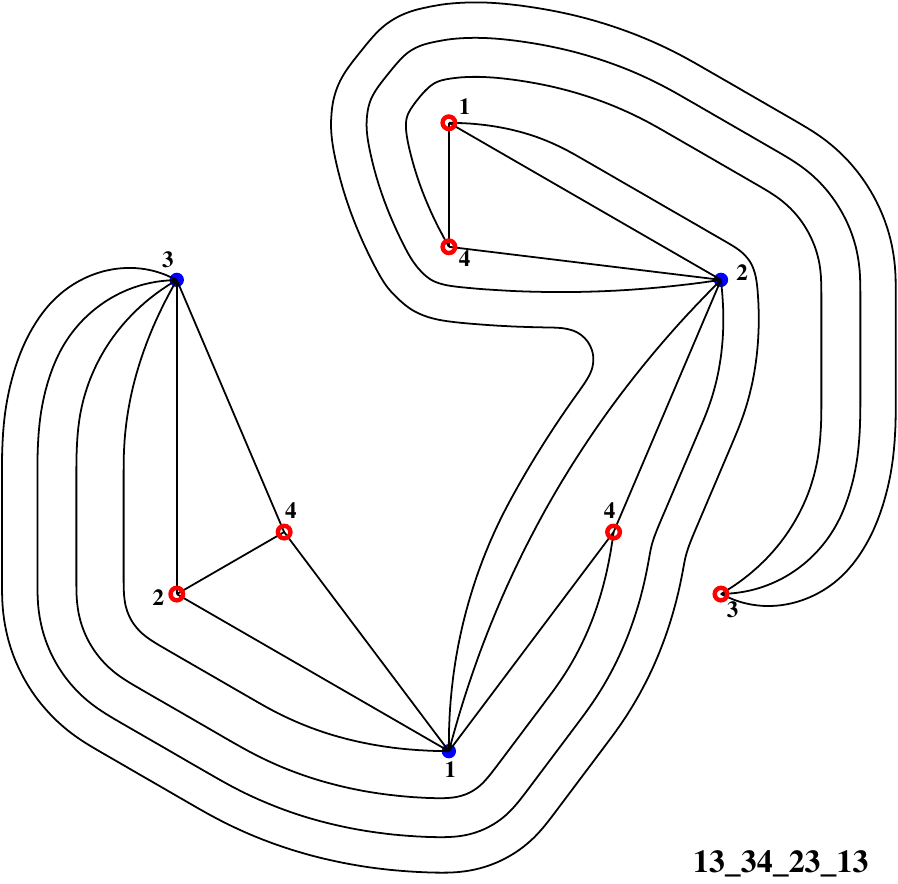}{(b) After $13,34,23,13$.}
 \caption{The last two moves of $P$; its first two stages are Figure~\ref{fig:d3-r1q1}(a), (c).}
 \label{fig:d3-p-stages}
\end{figure}

\begin{figure}[htbp]
 \centering
 \includegraphics[width=.50\textwidth]{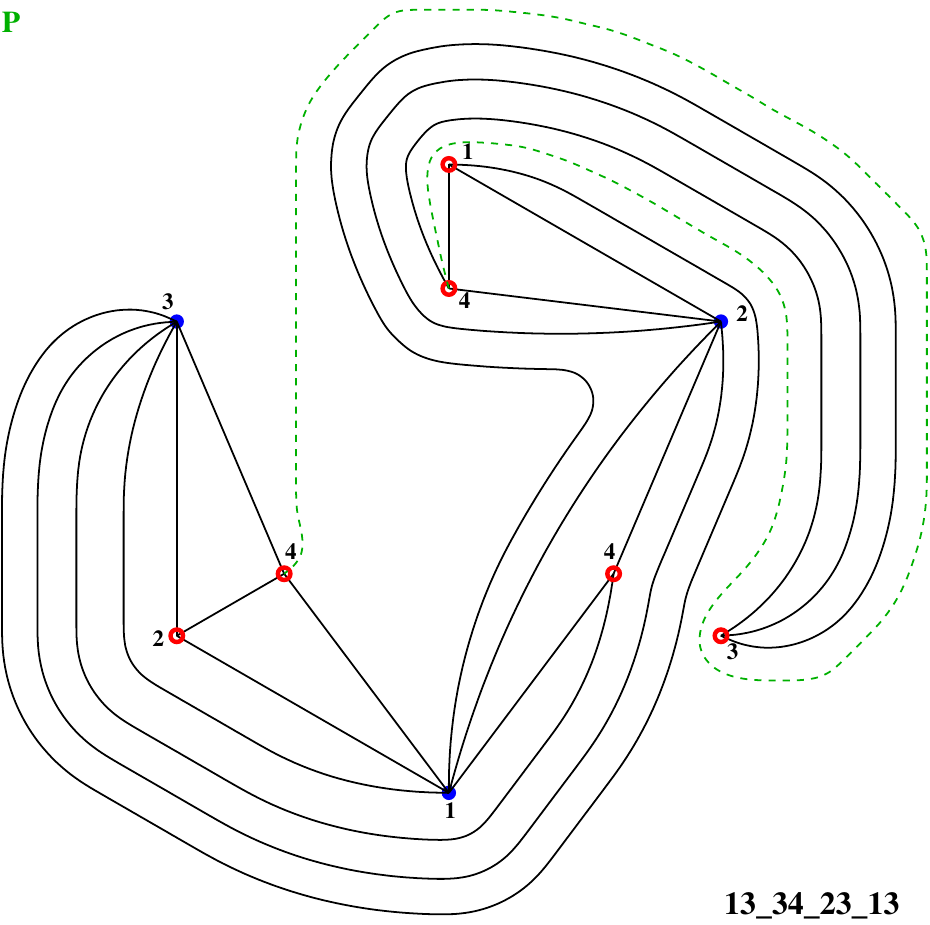}
 \caption{The terminal arc $\alpha_P$ from $a$ to $b$.}
 \label{fig:d3-p}
\end{figure}

For $P$, perform two further moves, $23,13$ (Figure~\ref{fig:d3-p-stages}). The terminal $4$-corners are $a,b$, joined by the arc $\alpha_P$ in Figure~\ref{fig:d3-p}. It has the same endpoints as $\alpha_{R_2}$ but winds differently, so $\delta_P$ and $\delta_{R_2}$ are different curves.

\FloatBarrier

\subsection{The seven vanishing cycles}\label{subsec:d3-cycles}
\begin{figure}[htbp]
 \centering
 \includegraphics[width=.38\textwidth]{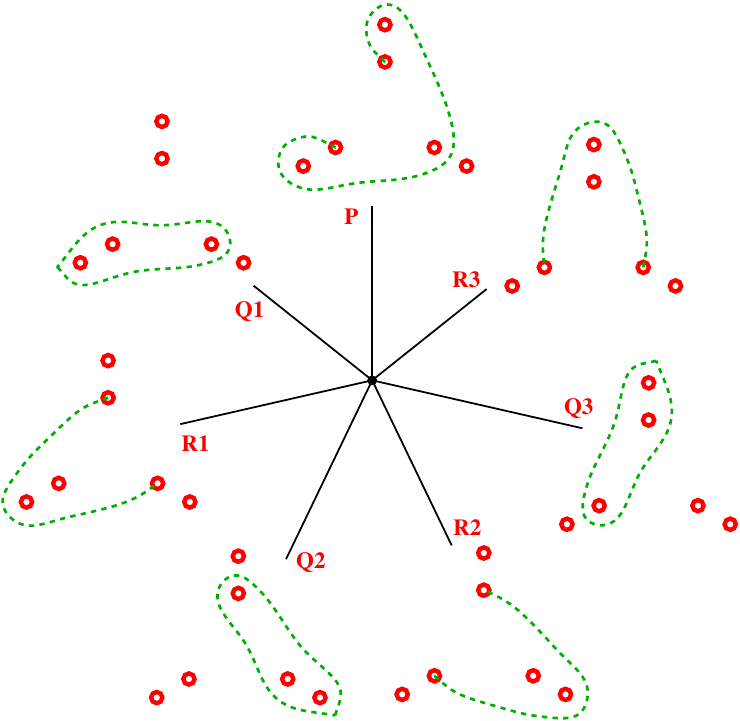}
 \caption{The seven degree-three vanishing curves. The genus $2$ surface is constructed by taking the double branched cover of the sphere branched over the six red points. The curve is the preimage of the green arc or circle.
   The arc give non-separating vanishing cycles while the circles give separating vanishing cycles.}
 \label{fig:d3-summary}
\end{figure}

Rotating $R_1,Q_1$ gives $R_2,R_3,Q_2,Q_3$. The lifts of the seven curves in Figure~\ref{fig:d3-summary} are the vanishing cycles: four nonseparating and three separating. The arc $\alpha_P$ is drawn differently in Figures~\ref{fig:d3-p} and~\ref{fig:d3-summary}; the two arcs are isotopic in $S$ relative to $W$. Their Dehn twists, in the order~\eqref{eq:d3order}, give the relation of Section~\ref{sec:relations}.

\section{The vanishing cycles in degree four}\label{sec:cycles}\label{subsec:d4}

The thirteen vanishing cycles come from five transports, $X_1,M_1,P_2,Y_2,M_0$, and the rotations of the first four. Unlike in degrees three and five, the base point $b_0$ is a vertex of the tiling.

\subsection{The reference cover and paths}\label{subsec:d4-reference}\label{subsec:d4-base}
\begin{figure}[htbp]
 \centering
 \includegraphics[width=.52\textwidth]{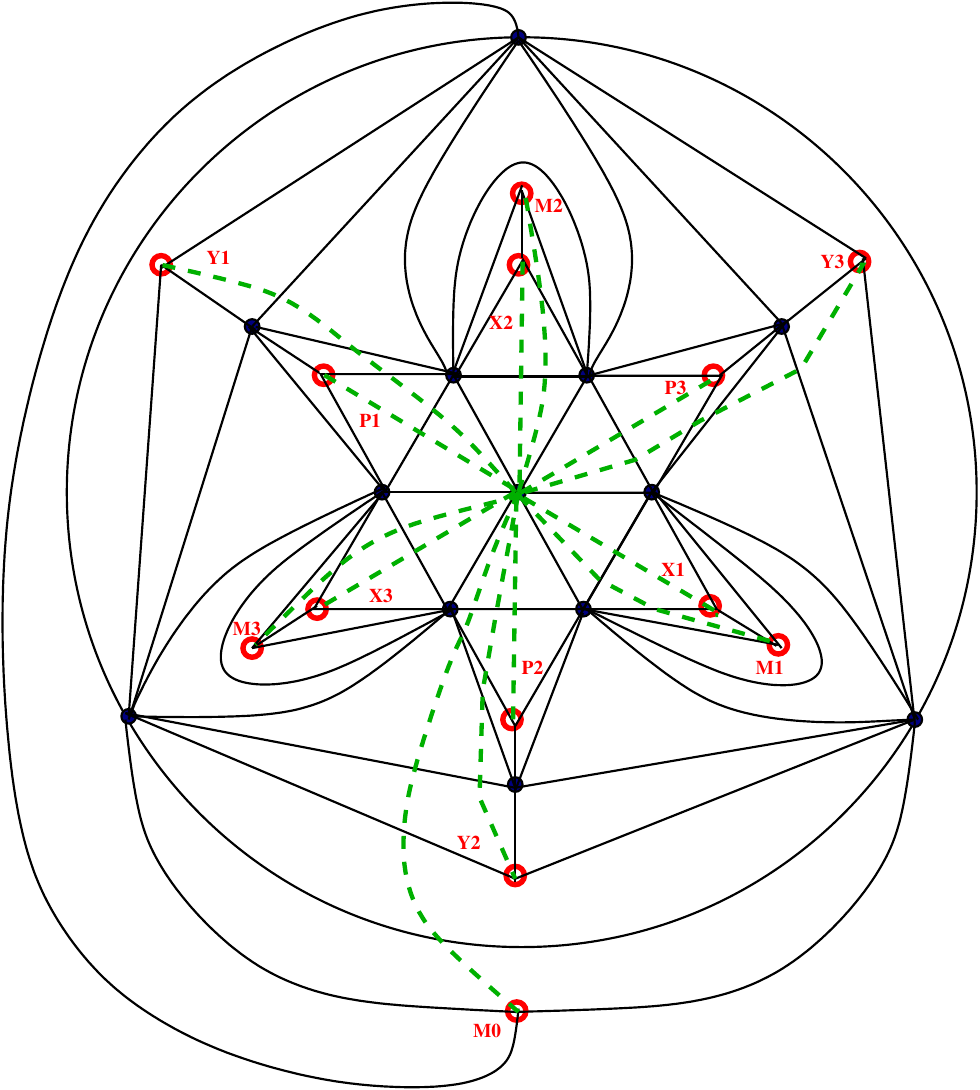}
 \caption{The thirteen degree-four vanishing paths, based at the central vertex $b_0$.}
 \label{fig:paths}
\end{figure}

At $b_0$ the two branch values coincide, so the reference cover has no exceptional hexagon (Section~\ref{sec:covers}). We therefore move from $b_0$ to a nearby point $b_0'$ over the face $234$, in the triangle through which $X_1$ and $M_1$ leave $b_0$ in Figure~\ref{fig:paths}; two triangles over $234$ meet at $b_0$, since $G$ has local degree two there. The resulting picture, marked by transport from $b_0$, is Figure~\ref{fig:first-move}(a).
The five representative paths in Figure~\ref{fig:paths} have the crossing words of Table~\ref{tab:d4-words}.

\begin{table}[!htbp]
\centering\small
\begin{tabular}{@{}cccc@{}}
\toprule
path & crossing word & terminal label & terminal curve\\
\midrule
$X_1$ & $23$ & $1$ & arc\\
$M_1$ & $23,13$ & $4$ & circle\\
$P_2$ & $34,13$ & $2$ & circle\\
$Y_2$ & $34,13,12,14,34$ & $2$ & arc\\
$M_0$ & $34,13,12,14,34,23,13$ & $4$ & circle\\
\bottomrule
\end{tabular}
\medskip
\caption{The five representative degree-four transports; the other words are obtained by rotation.}
\label{tab:d4-words}
\end{table}

For the rotated paths, rotate $b_0'$ as well and apply $1\mapsto2\mapsto3\mapsto1$ to the words. The subscript of a path increases by one (mod~$3$) under this rotation; for $X$, $P$ and $Y$ it is the terminal label. We take the thirteen paths in the order

\begin{equation}\label{eq:d4order}
 M_3,\ X_3,\ M_0,\ Y_2,\ P_2,\ M_1,\ X_1,\ Y_3,\ P_3,\ M_2,\ X_2,\ Y_1,\ P_1.
\end{equation}
Up to cyclic order, this is $M_0$ followed by the block $Y_2,P_2,M_1,X_1$ and its two rotations.

\subsection{The five transports}\label{subsec:d4-x1m1}\label{subsec:d4-p2}\label{subsec:d4-y2m0}
\begin{figure}[htbp]
 \centering
 \panel{.34}{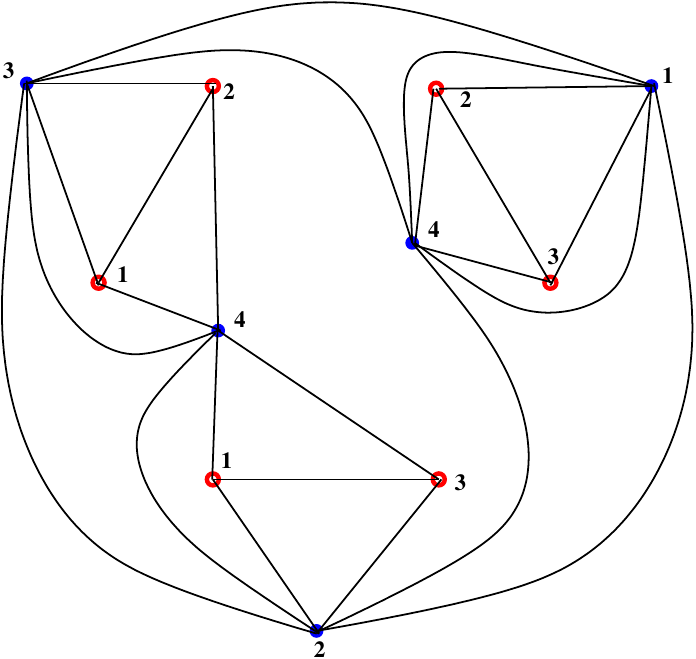}{(a) The perturbed starting picture, before crossing $23$.}\hspace{.08\textwidth}
 \panel{.34}{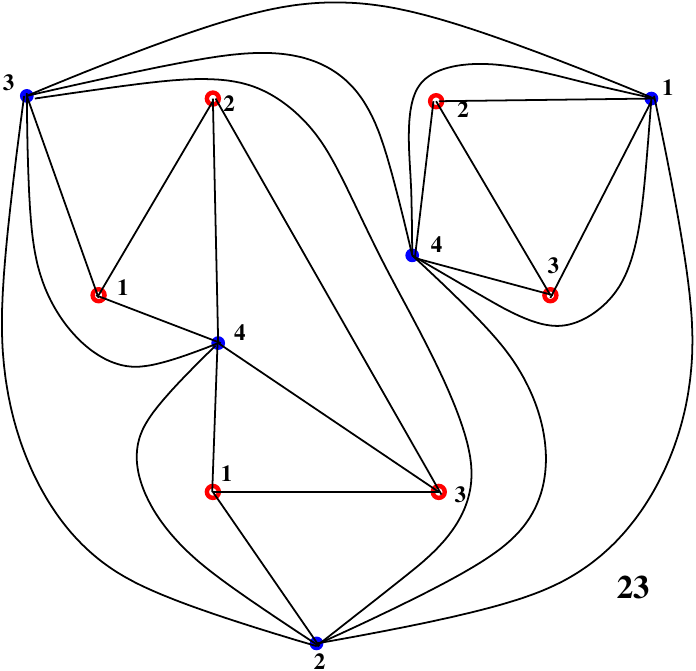}{(b) After crossing $23$.}
 \caption{The first move, $23$, from the perturbed reference picture.}
 \label{fig:first-move}
\end{figure}

\begin{figure}[htbp]
 \centering
 \panel{.34}{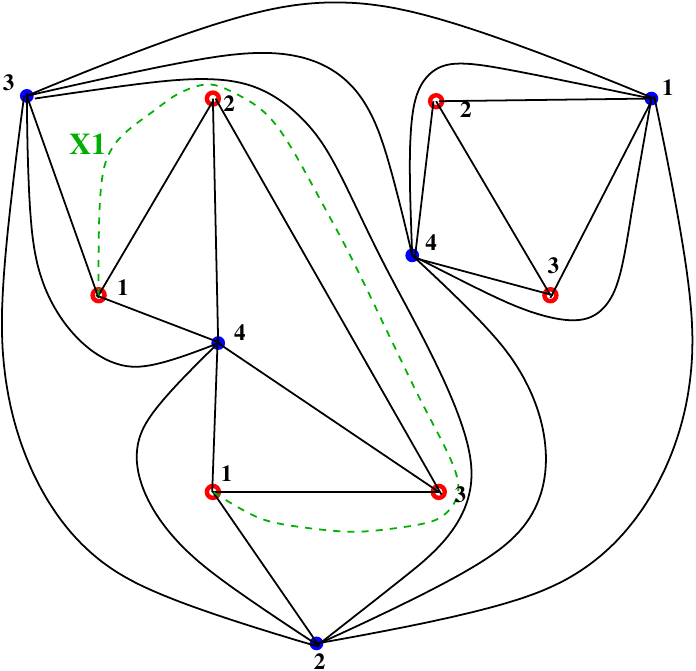}{(a) $X_1$: the arc from $r_1$ to $r'_1$.}\hspace{.08\textwidth}
 \panel{.34}{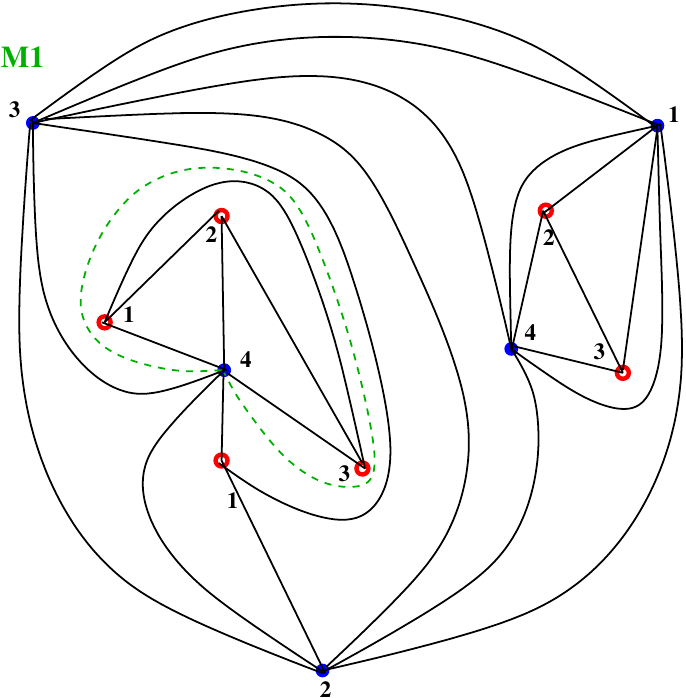}{(b) $M_1$: the circle through the left blue $4$.}
 \caption{The terminal curves of $X_1$ and $M_1$, after $23$ and $23,13$.}
 \label{fig:x1-m1}
\end{figure}

The first move, $23$, is shown in Figure~\ref{fig:first-move}. It gives the terminal arc $\alpha_{X_1}$ from $r_1$ to $r'_1$. A further move $13$ gives the circle $\gamma_{M_1}$ through the left blue $4$ (there are two blue vertices over $4$), separating $\{r_1,r'_2,r_3\}$ from $\{r'_1,r_2,r'_3\}$ (Figure~\ref{fig:x1-m1}).

\begin{figure}[htbp]
 \centering
 \includegraphics[width=.36\textwidth]{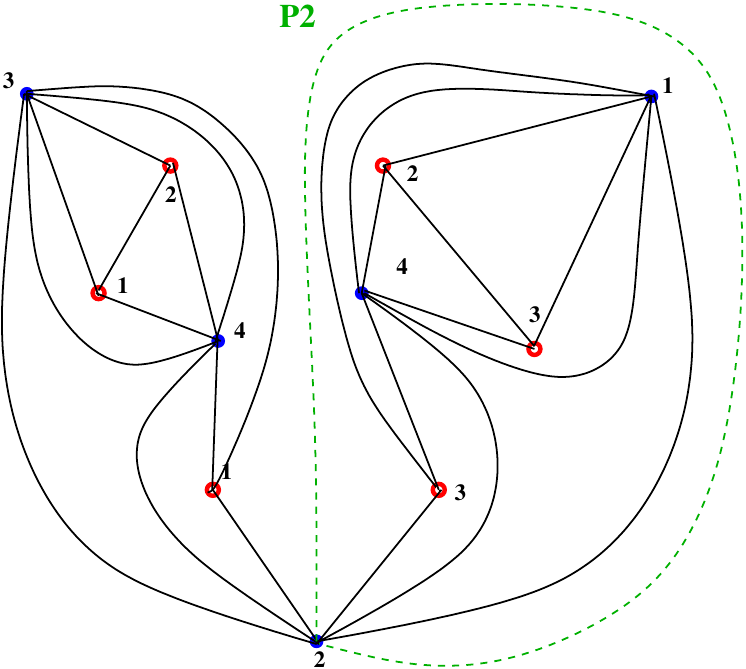}
 \caption{The terminal circle of $P_2$, after $34,13$.}
 \label{fig:p2}
\end{figure}

For $P_2$, the moves $34,13$ give the circle $\gamma_{P_2}$ through the blue $2$, separating $\{r_2,r_3,r'_3\}$ from $\{r_1,r'_1,r'_2\}$ (Figure~\ref{fig:p2}).

\begin{figure}[htbp]
 \centering
 \panel{.37}{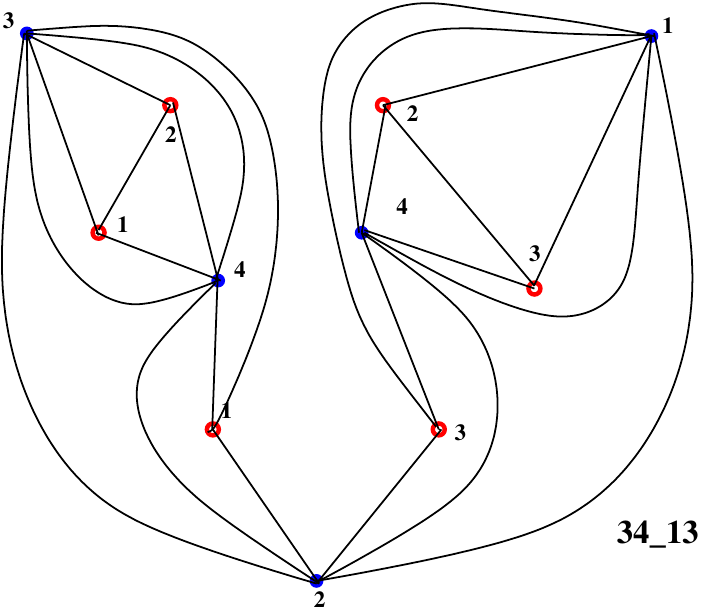}{(a) $34,13$.}\hspace{.08\textwidth}
 \panel{.37}{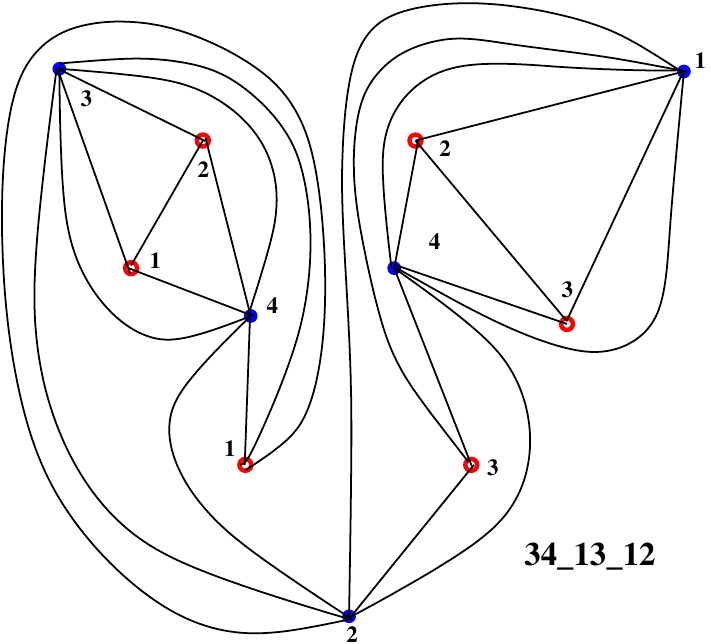}{(b) $34,13,12$.}
 \par\medskip
 \panel{.37}{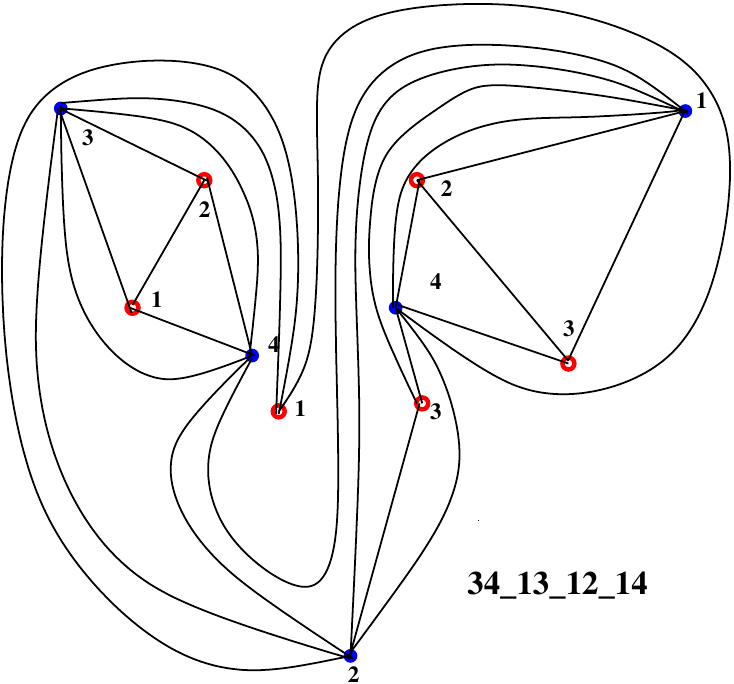}{(c) $34,13,12,14$.}\hspace{.08\textwidth}
 \panel{.37}{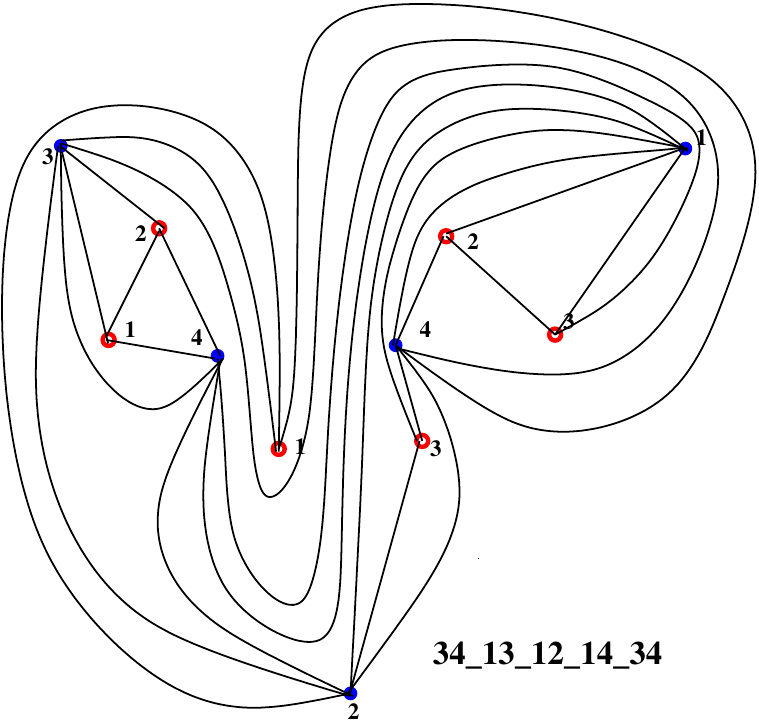}{(d) $34,13,12,14,34$.}
 \caption{The successive pictures along $Y_2$.}
 \label{fig:early-sequence}
\end{figure}

\begin{figure}[htbp]
 \centering
 \includegraphics[width=.38\textwidth]{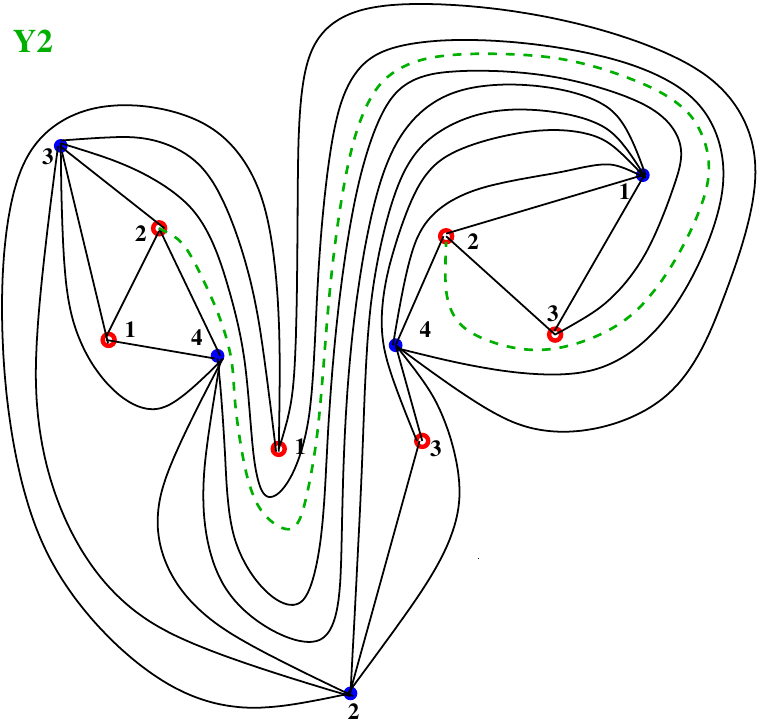}
 \caption{The terminal arc $\alpha_{Y_2}$ from $r'_2$ to $r_2$.}
 \label{fig:y2}
\end{figure}

Continuing from the picture of $P_2$ with $12,14,34$ gives $\alpha_{Y_2}$ from $r'_2$ to $r_2$ (Figures~\ref{fig:early-sequence}--\ref{fig:y2}). It has the same endpoints as $\alpha_{X_2}$ but winds differently, so $\delta_{Y_2}$ and $\delta_{X_2}$ are different curves.

\begin{figure}[htbp]
 \centering
 \panel{.40}{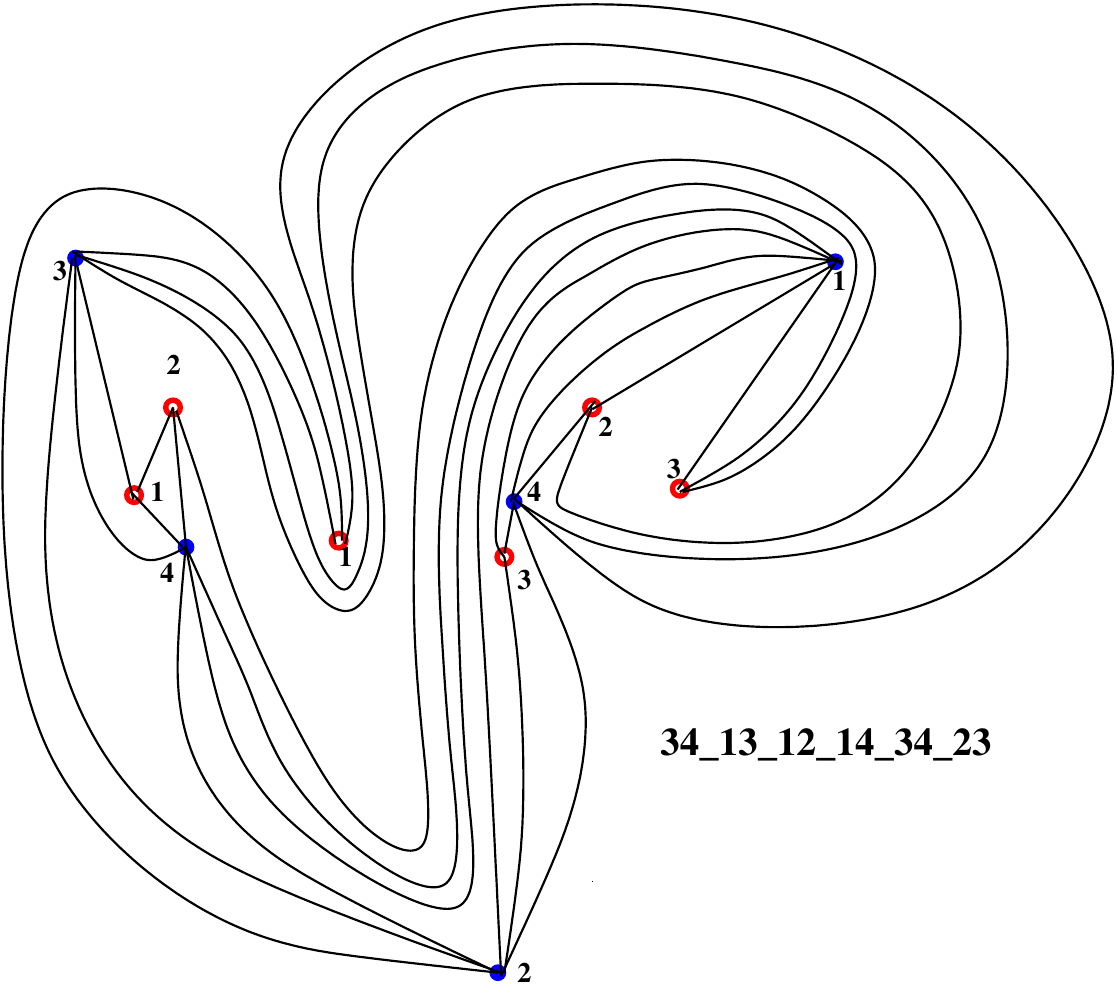}{(a) After the crossing $23$.}\hfill
 \panel{.40}{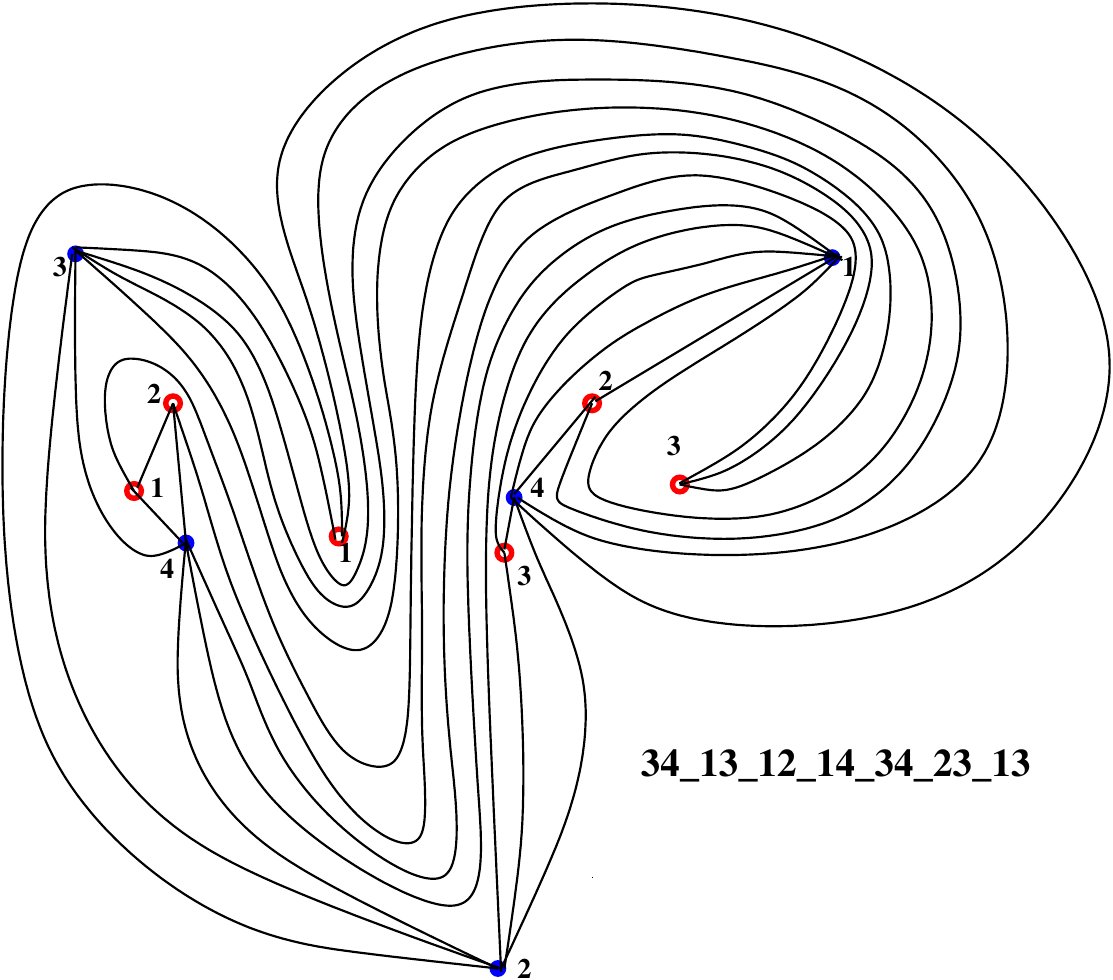}{(b) After the final crossing $13$.}
 \caption{The last two moves along $M_0$.}
 \label{fig:late-sequence}
\end{figure}

\begin{figure}[htbp]
 \centering
 \begin{minipage}[b]{.50\textwidth}\centering\setlength{\captionindent}{0pt}
  \includegraphics[width=\linewidth]{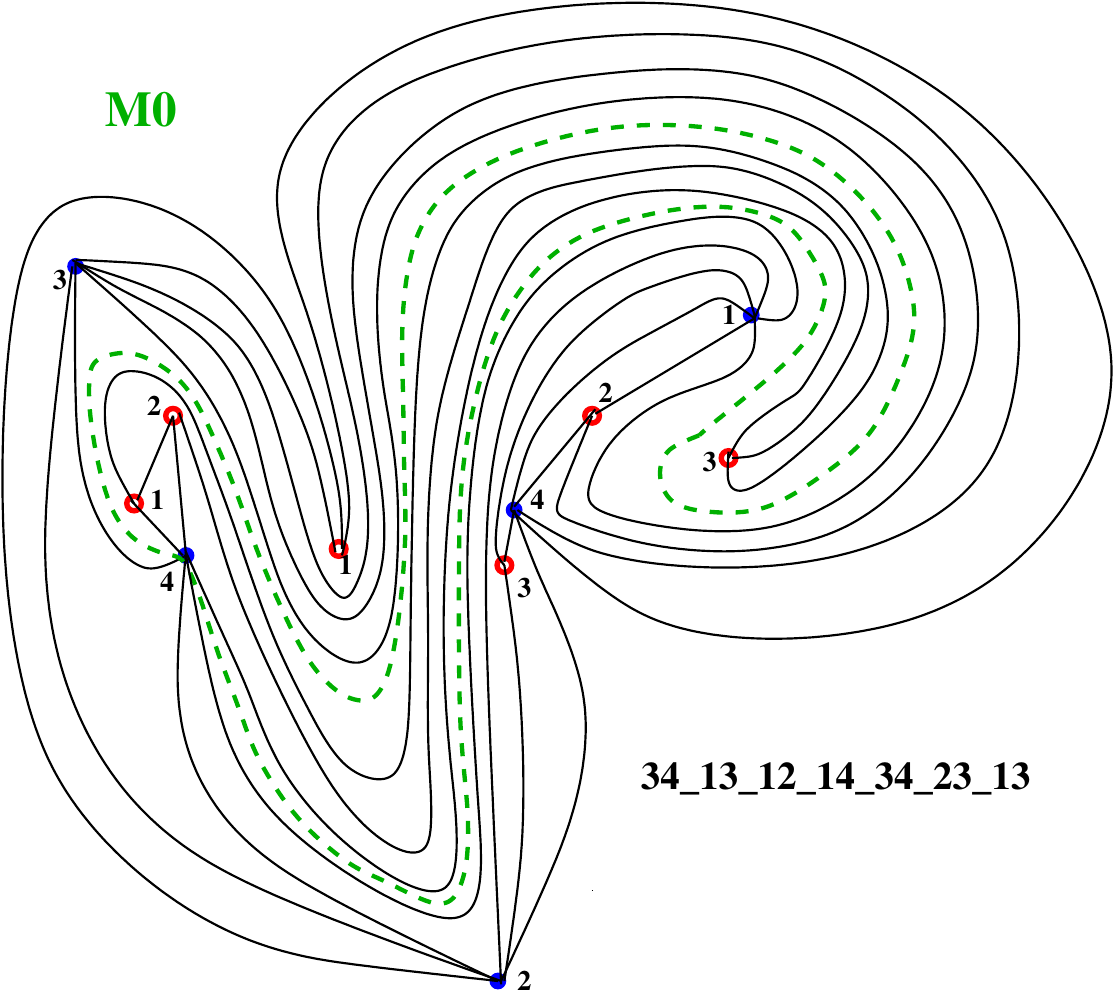}
  \caption{The terminal circle $\gamma_{M_0}$, after $34,13,12,14,34,23,13$.}
  \label{fig:m0}
 \end{minipage}\hfill
 \begin{minipage}[b]{.44\textwidth}\centering\setlength{\captionindent}{0pt}
  \includegraphics[width=\linewidth]{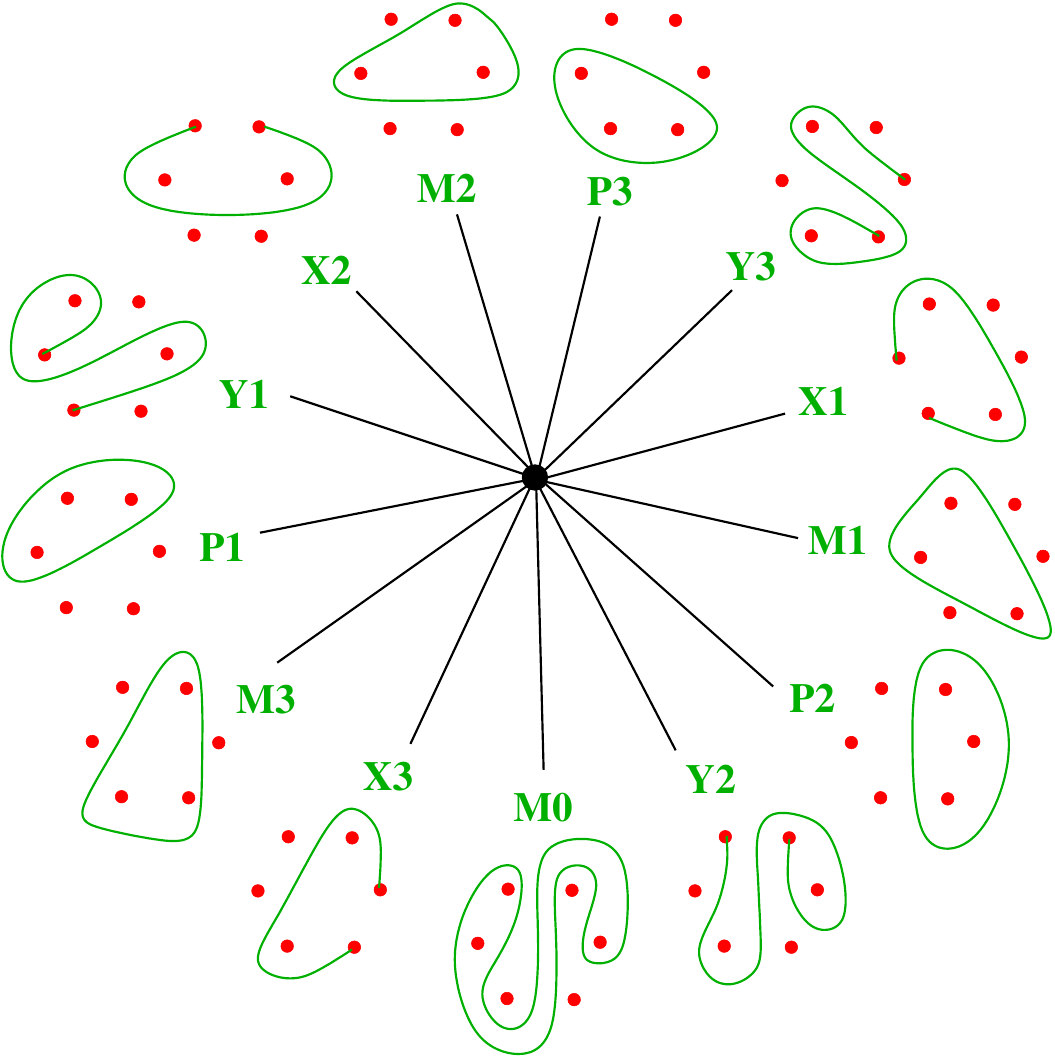}
  \caption{The thirteen degree-four vanishing curves.}
  \label{fig:all-cycles}
 \end{minipage}
\end{figure}

The last two moves, $23,13$, give $\gamma_{M_0}$ through the left blue $4$, separating $\{r_1,r'_2,r'_3\}$ from $\{r_2,r_3,r'_1\}$ (Figures~\ref{fig:late-sequence}--\ref{fig:m0}).

\subsection{The thirteen vanishing cycles}\label{subsec:d4-cycles}
Rotating the four representatives $X_1,Y_2,P_2,M_1$ gives the remaining eight curves. The lifts of the curves in Figure~\ref{fig:all-cycles} are the vanishing cycles: the arcs $\alpha_{X_i},\alpha_{Y_i}$ join $r_i$ to $r'_i$ and give six nonseparating cycles; the circles $\gamma_{P_i},\gamma_{M_i}$ ($i=1,2,3$) and $\gamma_{M_0}$ give seven separating cycles. Their Dehn twists, in the order~\eqref{eq:d4order}, give the relation of Section~\ref{sec:relations}.

\section{The vanishing cycles in degree five}\label{subsec:d5}

We apply the same construction in degree five, using a computer for the longer transports.

\subsection{The reference cover and paths}\label{subsec:d5-reference}\label{subsec:d5-base}
\begin{figure}[htbp]
 \centering
 \begin{minipage}[t]{.235\textwidth}\centering\includegraphics[width=0.878\linewidth]{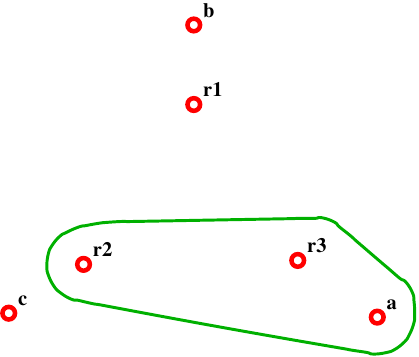}\par\smallskip{\small $A_1$}\end{minipage}\hspace{2.5em}
 \begin{minipage}[t]{.235\textwidth}\centering\includegraphics[width=0.838\linewidth]{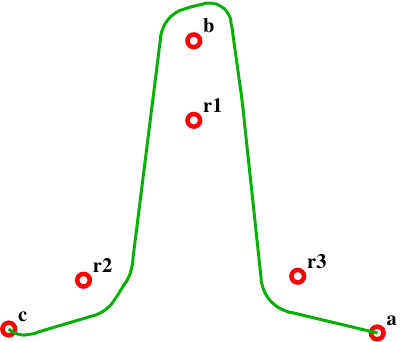}\par\smallskip{\small $A_2$}\end{minipage}\hspace{2.5em}
 \begin{minipage}[t]{.235\textwidth}\centering\includegraphics[width=1.000\linewidth]{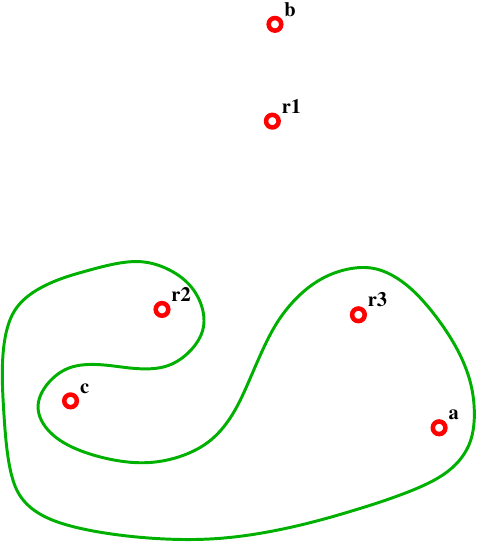}\par\smallskip{\small $A_3$}\end{minipage}\\[8pt]
 \begin{minipage}[t]{.235\textwidth}\centering\includegraphics[width=0.878\linewidth]{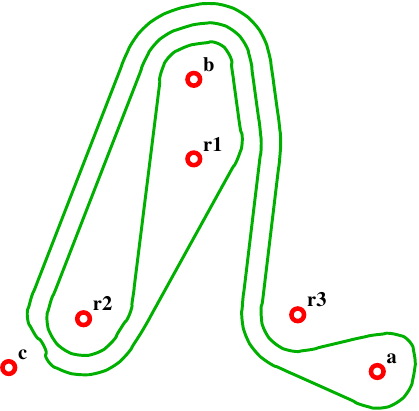}\par\smallskip{\small $A_4$}\end{minipage}\hspace{2.5em}
 \begin{minipage}[t]{.235\textwidth}\centering\includegraphics[width=0.917\linewidth]{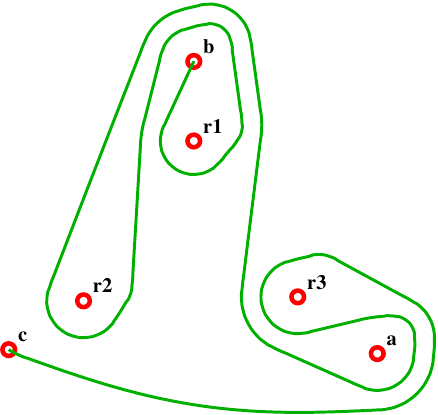}\par\smallskip{\small $A_5$}\end{minipage}\hspace{2.5em}
 \begin{minipage}[t]{.235\textwidth}\centering\includegraphics[width=0.878\linewidth]{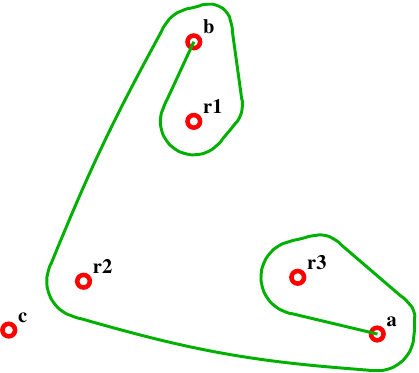}\par\smallskip{\small $A_6$}\end{minipage}\\[8pt]
 \begin{minipage}[t]{.235\textwidth}\centering\includegraphics[width=0.917\linewidth]{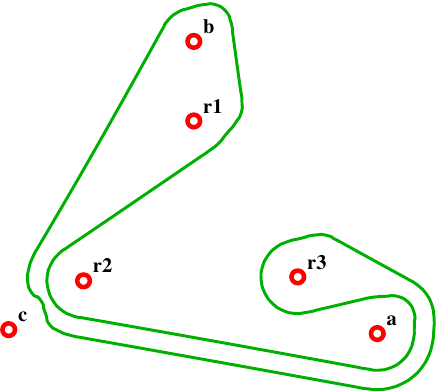}\par\smallskip{\small $A_7$}\end{minipage}\hspace{2.5em}
 \begin{minipage}[t]{.235\textwidth}\centering\includegraphics[width=0.838\linewidth]{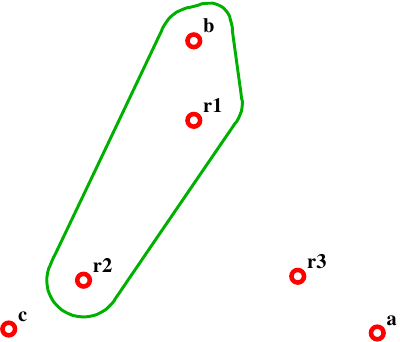}\par\smallskip{\small $A_{10}$}\end{minipage}\hspace{2.5em}
 \begin{minipage}[t]{.235\textwidth}\centering\includegraphics[width=0.917\linewidth]{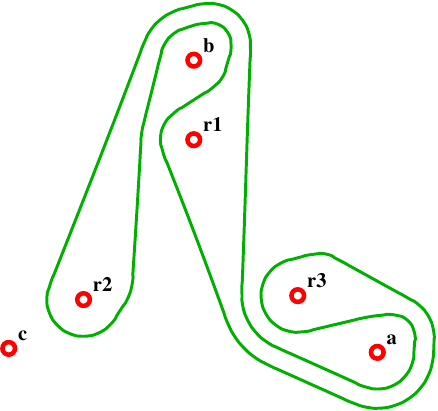}\par\smallskip{\small $A_{\mathrm{fix}}$}\end{minipage}\\[14pt]
 \begin{minipage}[t]{.48\textwidth}\centering\includegraphics[width=0.673\linewidth]{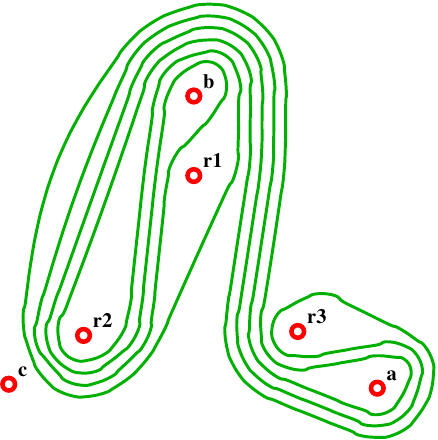}\par\smallskip{\small $A_8$}\end{minipage}\hfill
 \begin{minipage}[t]{.48\textwidth}\centering\includegraphics[width=0.654\linewidth]{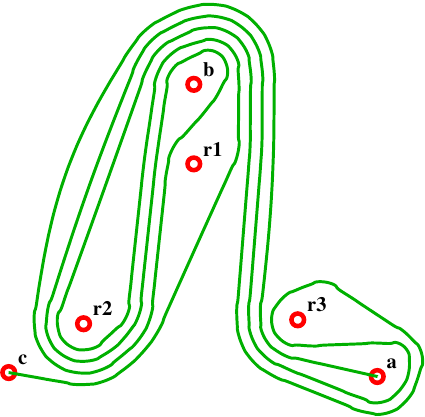}\par\smallskip{\small $A_9$}\end{minipage}
 \caption{The eleven representative vanishing curves, with the red points labelled. The remaining twenty curves are rotations of $A_1,\ldots,A_{10}$.}
 \label{fig:d5-summary}
\end{figure}

The reference picture is Figure~\ref{fig:intro-central}(c), with exceptional hexagon the outer face and base point $b_0$ at the center of the tiling. The red points $r_1,r_2,r_3$ are at the top, lower left, and lower right; the adjacent red points over $4$ are $b,c,a$, respectively, as labelled in Figure~\ref{fig:d5-summary}.

\begin{figure}[htbp]
 \centering
 \includegraphics[width=.5\textwidth]{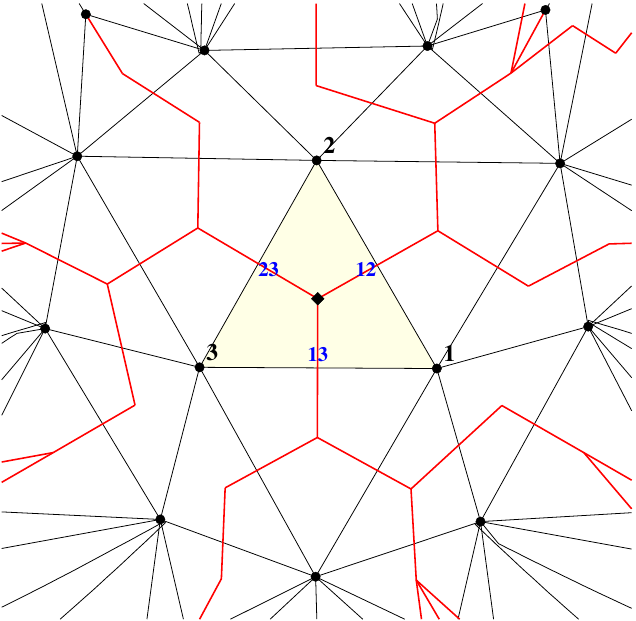}
 \caption{The labelled central triangle and the initial branches of the degree-five vanishing paths.}
 \label{fig:d5-centre}
\end{figure}

The thirty finite critical values form ten rotation orbits; the remaining value is at infinity. Choose the eleven representative paths $A_1,\ldots,A_{10},A_{\mathrm{fix}}$ with crossing words in Table~\ref{tab:d5-reps}. They follow shortest routes in the dual graph, starting from the labelled triangle of Figure~\ref{fig:d5-centre}. Write $A_k^{(j)}$ for the path obtained by replacing each label in the word $w_k$ of $A_k=A_k^{(0)}$ by its image under $\sigma^j$, where $\sigma\colon1\mapsto2\mapsto3\mapsto1$ fixes $4$. The resulting paths follow the rooted tree of Figure~\ref{fig:d5-paths}; separating their shared segments gives disjoint vanishing paths.

\subsection{Transport and the terminal curves}\label{subsec:d5-transport}\label{subsec:d5-pictures}
The program applies the hexagon move in the marked reference sphere. It records each edge's crossings with the reference decomposition and their order along each reference edge; this keeps track of how a new edge winds when it runs along the sides of the old hexagon. The terminal arc or circle is read as in Section~\ref{subsec:terminal}.

\begin{figure}[htbp]
 \centering
 \panel{.30}{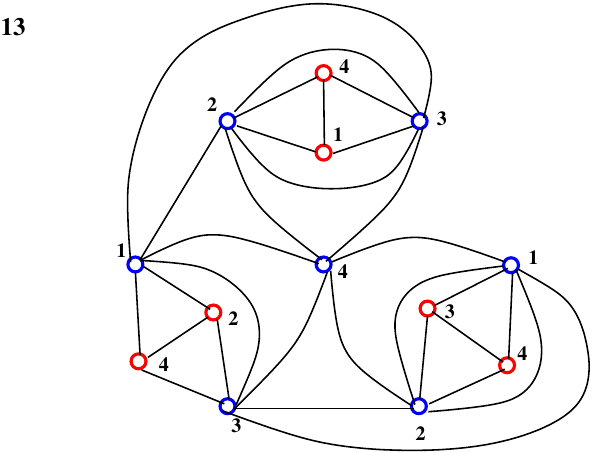}{$13$}\hfill
 \panel{.30}{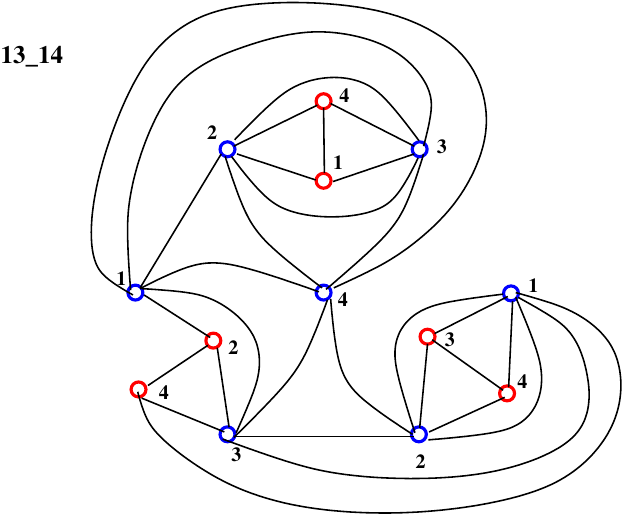}{$13,14$}\hfill
 \panel{.30}{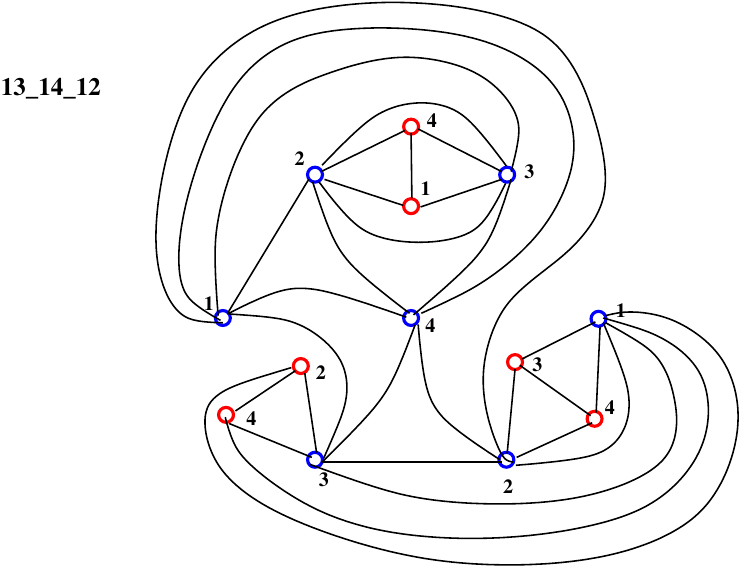}{$13,14,12$}
 \caption{The common initial stages of $A_{10},A_2,A_3$; the third stage belongs to $A_2,A_3$.}
 \label{fig:d5-first-route}
\end{figure}

\begin{figure}[htbp]
 \centering
 \panel{.31}{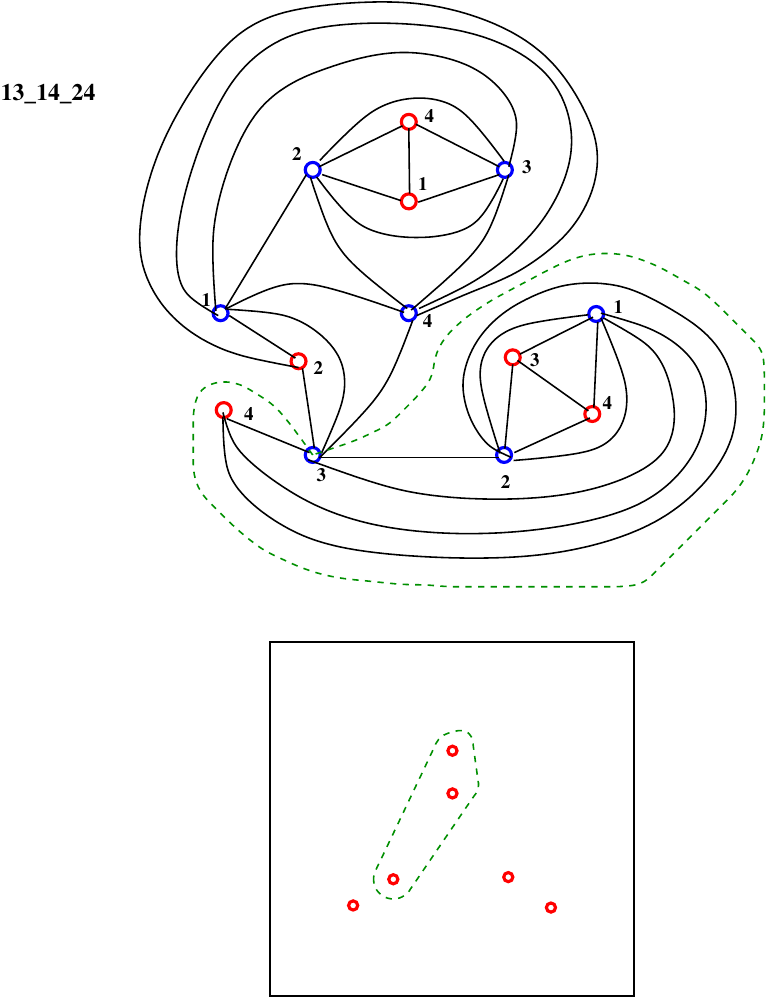}{(a) $A_{10}$: after $13,14,24$.}\hfill
 \panel{.31}{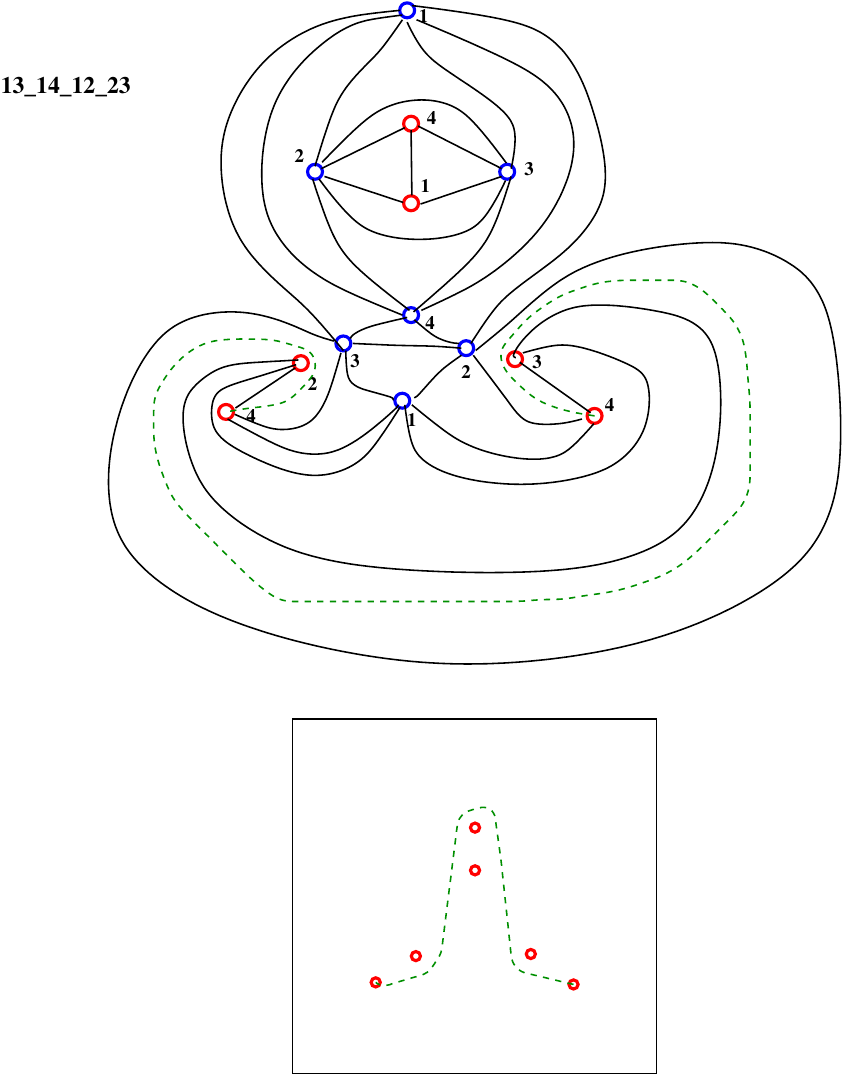}{(b) $A_2$: after $13,14,12,23$.}\hfill
 \panel{.31}{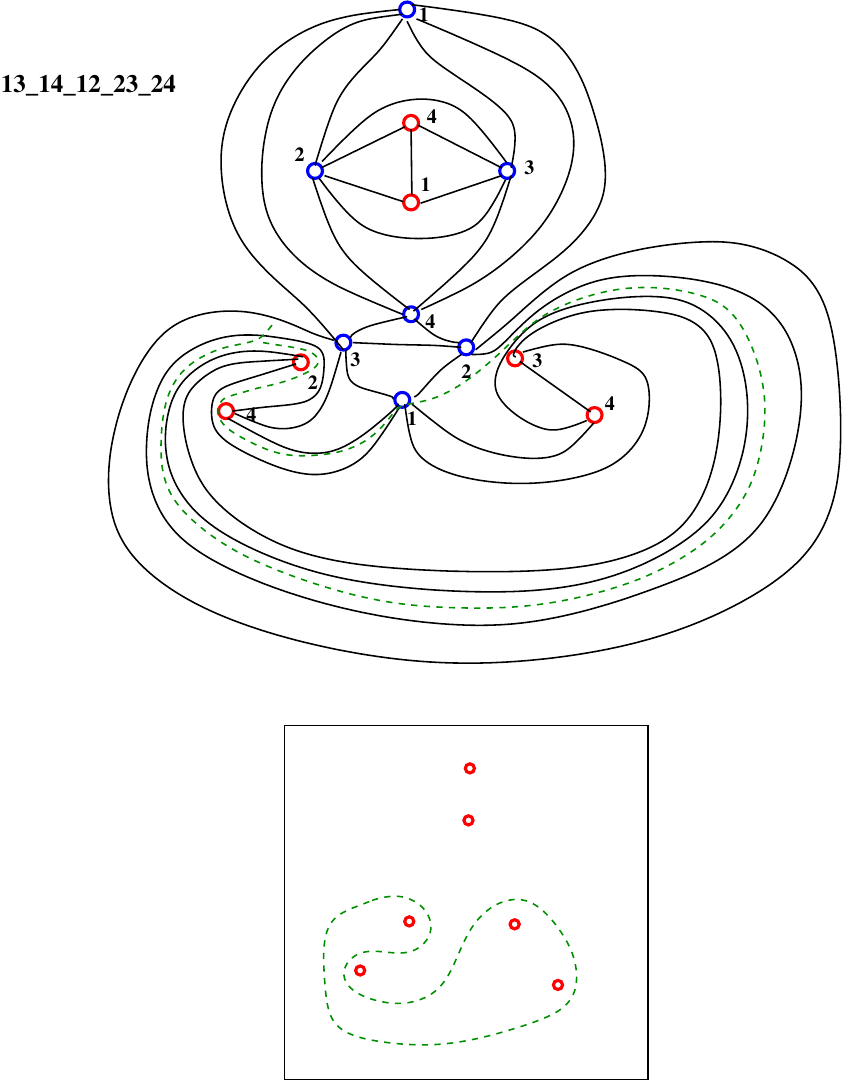}{(c) $A_3$: after $13,14,12,23,24$.}
 \caption{The terminal pictures of $A_{10},A_2,A_3$, with their arcs or circles below. Panel (b) is also the penultimate stage of $A_3$.}
 \label{fig:d5-sep-final}
\end{figure}

Figures~\ref{fig:d5-first-route} and~\ref{fig:d5-sep-final} give three complete transports. After the common moves $13,14$, crossing $24$ gives $A_{10}$: the two $3$-corners coincide at a blue vertex, and the terminal circle encloses $\{b,r_1,r_2\}$. Instead crossing $12,23$ gives $A_2$, whose terminal arc joins the red $4$-corners $c,a$. One further crossing, $24$, gives $A_3$: the two $1$-corners coincide, and the terminal circle encloses $\{a,r_2,r_3\}$.

\FloatBarrier

\begin{table}[!htbp]
\centering\small
\setlength{\tabcolsep}{4pt}
\begin{tabular}{@{}clclc@{}}
\toprule
 & crossing word $w_k$ & label & terminal curve & type\\
\midrule
$A_1$ & $13,34,24$ & $1$ & circle: $\{a,r_2,r_3\}\mid\{b,c,r_1\}$ & separating\\
$A_2$ & $13,14,12,23$ & $4$ & arc $c$--$a$ & nonseparating\\
$A_3$ & $13,14,12,23,24$ & $1$ & circle: $\{a,r_2,r_3\}\mid\{b,c,r_1\}$ & separating\\
$A_4$ & $13,14,24,23,12,14,34,24$ & $1$ & circle: $\{c,r_2,r_3\}\mid\{a,b,r_1\}$ & separating\\
$A_5$ & $13,14,24,23,12,14,34,23,12$ & $4$ & arc $b$--$c$ & nonseparating\\
$A_6$ & $13,14,24,34,14,12,23$ & $4$ & arc $a$--$b$ & nonseparating\\
$A_7$ & $13,14,24,34,14,12,23,34,13$ & $2$ & circle: $\{a,c,r_2\}\mid\{b,r_1,r_3\}$ & separating\\
$A_8$ & $13,14,24,23,12,14,34,23,13,14,24,34,13$ & $2$ & circle: $\{a,r_1,r_3\}\mid\{b,c,r_2\}$ & separating\\
$A_9$ & $13,14,24,23,12,14,34,23,13,14,24,23,13$ & $4$ & arc $a$--$c$ & nonseparating\\
$A_{10}$ & $13,14,24$ & $3$ & circle: $\{b,r_1,r_2\}\mid\{a,c,r_3\}$ & separating\\
$A_{\mathrm{fix}}$ & $13,14,24,23,12,14,34,23,13,14,12,23$ & $4$ & circle: $\{a,b,c\}\mid\{r_1,r_2,r_3\}$ & separating\\
\bottomrule
\end{tabular}
\medskip
\caption{The eleven representative transports. The remaining words and curves are obtained by rotation.}
\label{tab:d5-reps}
\end{table}
\FloatBarrier

\FloatBarrier

\subsection{The vanishing cycles and their order}\label{subsec:d5-reps}\label{subsec:d5-cycles}\label{subsec:d5-paths}
Rotate the ten curves $A_1,\ldots,A_{10}$ of Figure~\ref{fig:d5-summary} and adjoin $A_{\mathrm{fix}}$. The arcs $A_2,A_5,A_6,A_9$ give twelve nonseparating lifts; the remaining circles give nineteen separating lifts.

\begin{figure}[htbp]
 \centering
 \includegraphics[width=\textwidth,height=.55\textheight,keepaspectratio]{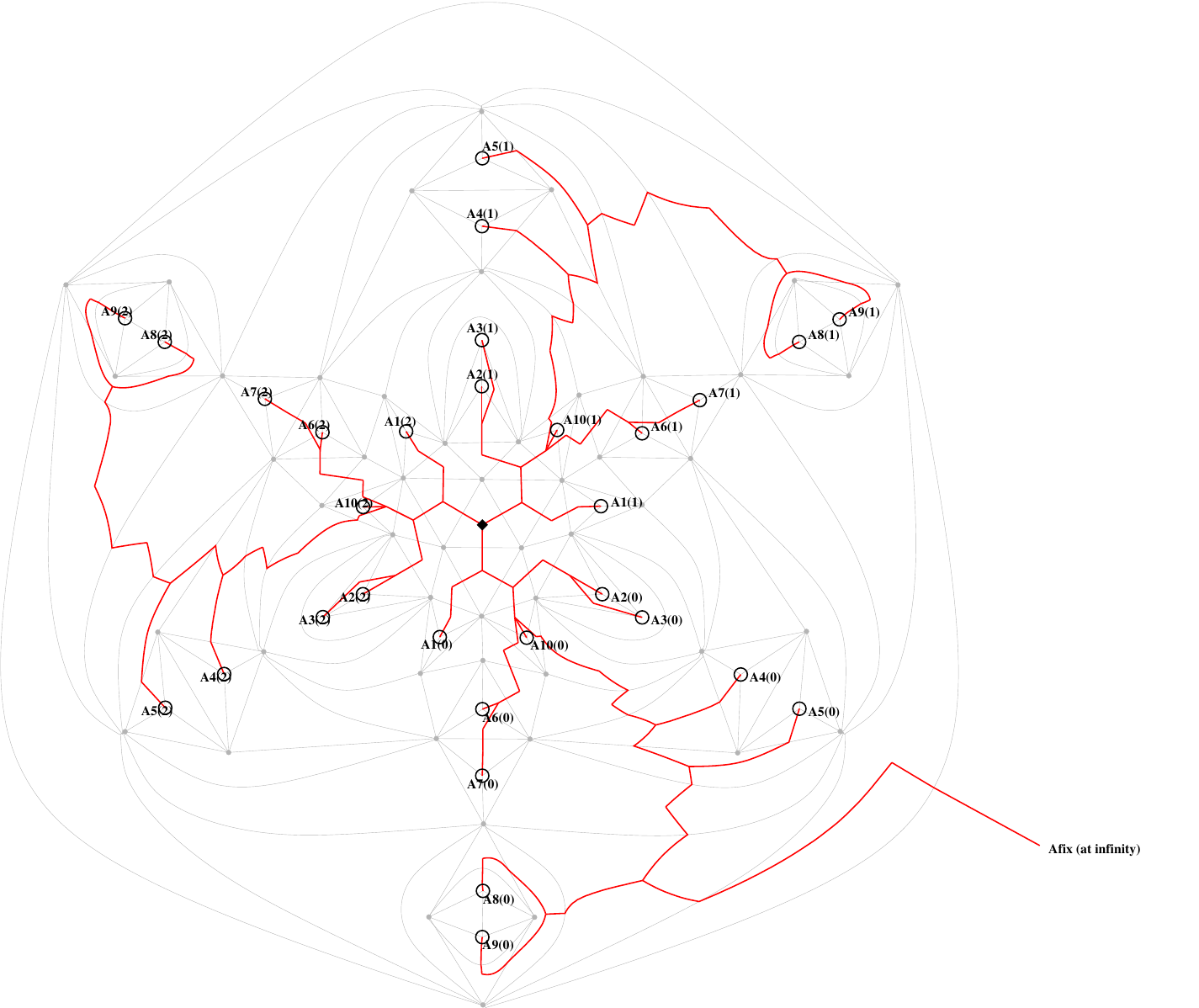}
 \caption{The thirty-one paths. The diamond is $b_0$; $A_{\mathrm{fix}}$ runs to infinity. Shared initial segments are drawn once.}
 \label{fig:d5-paths}
\end{figure}

We take the paths of Figure~\ref{fig:d5-paths} in the order
\begin{equation}\label{eq:d5order}
 A_1^{(j)},\ A_6^{(j)},\ A_7^{(j)},\ A_{10}^{(j)},\ A_8^{(j)},\ A_9^{(j)},\ A_5^{(j)},\ A_4^{(j)},\ A_3^{(j)},\ A_2^{(j)}\qquad(j=0,1,2\text{ in turn}),
\end{equation}
with $A_{\mathrm{fix}}$ inserted between $A_9^{(0)}$ and $A_5^{(0)}$. The three blocks exit through $13,12,23$, with the orientation of Figure~\ref{fig:d5-centre}.

\section{The relations}\label{sec:relations}

\subsection{The relations}\label{subsec:birman-hilden}\label{subsec:three-relations}\label{subsec:proof-B}
Theorem~B follows from the previous sections; no further computation is needed. Sections~\ref{subsec:d3}, \ref{subsec:d4-cycles} and~\ref{subsec:d5-cycles} identify the vanishing cycle $\delta_i$ along each vanishing path $\gamma_i$ of Figures~\ref{fig:d3-paths}, \ref{fig:paths} and~\ref{fig:d5-paths} as a curve in the fiber $\Sigma=f^{-1}(b_0)$. The monodromy along the loop that runs out along $\gamma_i$, once around its critical value in the positive sense, and back is the right-handed Dehn twist $t_{\delta_i}$. In the orders~\eqref{eq:d3order}, \eqref{eq:d4order} and~\eqref{eq:d5order} these loops compose to a loop around all the critical values. That loop bounds a disk in $B\cong\PP^1$ containing no critical value, so its monodromy is trivial: $t_{\delta_n}\circ\cdots\circ t_{\delta_1}=1$ in $\Mod(\Sigma)$.

Under the hyperelliptic quotient $\pi\colon\Sigma\to S$, the twist $t_\delta$ covers the half-twist $H_\alpha$ when $\delta=\pi^{-1}(\alpha)$ is nonseparating, and the squared twist $T_\gamma^2$ when $\delta=\pi^{-1}(\gamma)$ is separating~\cite{BH,FM}. So each factorization also gives an identity among half-twists and squared twists in $\Mod(S,W)$, the mapping class group of the sphere with the six red points marked.

\subsection{A direct check}\label{subsec:braids}\label{subsec:mcg-d3}\label{subsec:mcg-d4}\label{subsec:mcg-d5}
The relations can also be checked directly, in two steps. First, the half-twists and squared twists compose to the identity in $\Mod(S,W)$, the mapping class group of the sphere with the six red points marked. This is an identity in a quotient of the braid group $B_6$, and it can be checked by composing automorphisms of the free group $\pi_1(S\setminus W)$. By the theorem of Birman and Hilden~\cite{BH}, the product $t_{\delta_n}\circ\cdots\circ t_{\delta_1}$ is then $1$ or the hyperelliptic involution $\iota$. Second, the product acts trivially on $H_1(\Sigma;\ZZ)$, while $\iota$ acts as $-1$.

\section{The complex surface and Xiao's fibration}\label{sec:surface}

\subsection{The complex surface}\label{subsec:complex-base}\label{subsec:complex-surface}
Give $B$ the complex structure making $G\colon B\to T$ holomorphic~\cite{Forster}. In the hexagonal model $E\colon y^2=x^3-1$, with $q=x$, choose $t$ so that $b_0=0$ and the other fixed point of $\sigma$ is at infinity. Then $G(\omega t)=\omega G(t)$.

The curve $\Delta=\{G(t)=q(e)\}\subset B\times E$ is algebraic. Its local equation~\eqref{eq:local-branch} gives $0,4,24$ nodes, $3,9,27$ cusps, and $7,13,31$ simple tangencies to the fibers of $B\times E\to B$, respectively.

The reference covers  have
no nontrivial deck transformations. The identifications used
to glue the local families are therefore unique.

Away from $\Pi^{-1}(\Delta)$, pull back the complex structure of
$B\times E$. The simple-branch charts and the local models 
discussed earlier extend it over $\Pi^{-1}(\Delta)$. The transition maps
are holomorphic in the complement of this set and continuous across it, hence
holomorphic by removable singularities. Thus $X$ has a unique
complex structure making $\Pi$ holomorphic.

\subsection{Identification with Xiao's family}\label{subsec:identification}
The covers $h_t$ are primitive. This is automatic for $d=3,5$; for $d=4$, the Weierstrass distribution $(2,2,2,0)$ rules out a factorization through a degree-two isogeny, since such an isogeny maps $E'[2]$ to only two points. Xiao's universal family therefore gives a classifying map $\varphi\colon B\to X(d)$~\cite{Xiao,KaniHurwitz}.

The calculation of Section~\ref{subsec:base} identifies the generic points of $B$ with the generic covers in $X(d)$, and the local models extend this identification over the vertices. Hence $\varphi$ is an isomorphism and $f$ is Xiao's fibration, proving Theorem~A.

Alternatively, the degree of $\varphi$ follows from the singular-fiber count, without the exhaustive enumeration. Let $m=\deg\varphi$ and $\nu=7,13,31$. After resolving singularities, the pullback of Xiao's fibration by $\varphi$ has $m\nu$ nodes in its singular fibers. Thus
\[
 -4+\nu=e(X)=-4+m\nu,
\]
so $m=1$.

\subsection{The second projection}\label{subsec:invariants-surface}
In the Stein factorization $X\to E'\to E$ of $g$, the curve $E'$ is elliptic because $f$ is non-isotrivial, and primitivity forces $E'\to E$ to have degree one. Thus the generic fiber $F_e=g^{-1}(e)$ is connected. Its degree-$d$ map to $B$ has $\deg G$ simple branch points, giving
\[
 \operatorname{genus}(F_e)=\tfrac12\deg G-d+1=0,3,16.
\]
For $d=3$, $X$ is an elliptic ruled surface blown up three times, hence diffeomorphic to $(S^2\times T^2)\#3\overline{\mathbb{CP}}{}^2$~\cite{BK}. For $d=4$ we have $p_g=q=1$ and $K^2=3$, and $g$ is the Albanese fibration. These surfaces with a genus-two pencil were constructed by Xiao~\cite[Theorem~6.5]{Xiao}, and Polizzi~\cite{Pol} describes them as divisors in the symmetric cube of~$E$.

\subsection{Fundamental groups}\label{subsec:pi1}
For $d=3$, $X$ is an elliptic ruled surface blown up three times, so $\pi_1(X)\cong\ZZ^2$. For $d=4,5$, a direct computation from the vanishing cycles of Theorem~B gives $\pi_1(X)\cong\ZZ^2$ in both cases. In degree four this was first computed by \.{I}nan\c{c} Baykur (personal communication).

\subsection{Sections}\label{subsec:sections}
For $d=3,5$ the Weierstrass distribution is $(1,1,1,3)$, so each smooth fiber of $f$ has exactly one Weierstrass point over each of $e_1,e_2,e_3$; in the reference fiber these are the red points $r_1,r_2,r_3$. The hyperelliptic involutions of the fibers form an involution of $X$, and these points trace three components of its fixed curve, each mapping isomorphically to $B$. Thus the three red points give three holomorphic sections of $f$. They are disjoint, since $g$ maps them to the three different points $e_1,e_2,e_3$. In the case $d=4$ it is easy to see, using the projection $g$, that $f$ has no holomorphic section. Nevertheless, $f$ admits a smooth section. This can be checked by lifting the relation of Theorem~B to the mapping class group of $\Sigma$ with one marked point; we leave the verification to the interested reader.

\subsection{Other elliptic curves}\label{subsec:other-E}
The same construction applies to any elliptic curve $E$: only the complex structures on $B$ and $X$ change. The marked topological covering and its vanishing cycles remain unchanged, so Theorem~B and the fundamental groups above hold for every $E$. The hexagonal curve was chosen only for its order-three symmetry.


\begin{thebibliography}{99}
\bibitem{Akh} A.~Akhmedov, \emph{Xiao's genus-two fibration: branched covers and braid monodromy}, preprint, arXiv:2609.17634 (2026).
\bibitem{AM} A.~Akhmedov and N.~Monden, \emph{Genus $2$ Lefschetz fibrations with $b_2^+=1$ and $c_1^2=1,2$}, Kyoto J. Math. \textbf{60} (2020), 1419--1451.
\bibitem{Aur} D.~Auroux, \emph{Fiber sums of genus $2$ Lefschetz fibrations}, Turkish J. Math. \textbf{27} (2003), 1--10. 
\bibitem{BK} R.~\.I.~Baykur and M.~Korkmaz, \emph{Small Lefschetz fibrations and exotic $4$-manifolds}, Math. Ann. \textbf{367} (2017), 1333--1361.
\bibitem{BH} J.~S.~Birman and H.~M.~Hilden, \emph{On isotopies of homeomorphisms of Riemann surfaces}, Ann. of Math. (2) \textbf{97} (1973), 424--439.
\bibitem{En} H.~Endo, \emph{Meyer's signature cocycle and hyperelliptic fibrations}, Math. Ann. \textbf{316} (2000), 237--257.
\bibitem{FM} B.~Farb and D.~Margalit, \emph{A primer on mapping class groups}, Princeton Mathematical Series 49, Princeton University Press, 2012.
\bibitem{Forster} O.~Forster, \emph{Lectures on Riemann surfaces}, Graduate Texts in Mathematics 81, Springer-Verlag, New York, 1981. 
\bibitem{Fox} R.~H.~Fox, \emph{Covering spaces with singularities}, in: Algebraic geometry and topology. A symposium in honor of S.~Lefschetz, Princeton University Press, 1957, 243--257.
\bibitem{FK} G.~Frey and E.~Kani, \emph{Curves of genus $2$ covering elliptic curves and an arithmetical application}, in: Arithmetic algebraic geometry (Texel, 1989), Progr. Math. 89, Birkh\"auser, Boston, 1991, 153--176. 
\bibitem{GS} R.~E.~Gompf and A.~I.~Stipsicz, \emph{$4$-Manifolds and Kirby Calculus}, Graduate Studies in Mathematics 20, American Mathematical Society, 1999.
\bibitem{Huang} E.~Huang, \emph{Equivalent genus-$2$ factorizations of type $(4,3)$}, preprint, arXiv:2602.20451 (2026).
\bibitem{Kani} E.~Kani, \emph{The number of curves of genus two with elliptic differentials}, J. Reine Angew. Math. \textbf{485} (1997), 93--121. 
\bibitem{KaniHurwitz} E.~Kani, \emph{Hurwitz spaces of genus $2$ covers of an elliptic curve}, Collect. Math. \textbf{54} (2003), no.~1, 1--51.
\bibitem{Kas} A.~Kas, \emph{On the handlebody decomposition associated to a Lefschetz fibration}, Pacific J. Math. \textbf{89} (1980), 89--104.

\bibitem{Kaya}
G.~Karado\u{g}an-Kaya,
\emph{On the moduli of surfaces admitting genus two fibrations
over elliptic curves},
Arch. Math. \textbf{89} (2007), no.~4, 315--325.

\bibitem{Kor} M.~Korkmaz, \emph{Noncomplex smooth $4$-manifolds with Lefschetz fibrations}, Internat. Math. Res. Notices 2001, no.~3, 115--128.
\bibitem{Mat} Y.~Matsumoto, \emph{Lefschetz fibrations of genus two---a topological approach}, in: Topology and Teichm\"uller spaces (Katinkulta, 1995), World Scientific, 1996, 123--148.
\bibitem{Nak} K.~Nakamura, \emph{Geography of genus $2$ Lefschetz fibrations}, preprint, arXiv:1811.03708 (2018). 
\bibitem{OS} B.~Ozbagci and A.~I.~Stipsicz, \emph{Noncomplex smooth $4$-manifolds with genus-$2$ Lefschetz fibrations}, Proc. Amer. Math. Soc. \textbf{128} (2000), 3125--3128.
\bibitem{Pol} F.~Polizzi, \emph{On surfaces of general type with $p_g=q=1$, $K^2=3$}, Collect. Math. \textbf{56} (2005), no.~2, 181--234.
\bibitem{ST} B.~Siebert and G.~Tian, \emph{On the holomorphicity of genus two Lefschetz fibrations}, Ann. of Math. (2) \textbf{161} (2005), 959--1020.
\bibitem{Trieste} A.~I.~Stipsicz and Z.~Szab\'o, \emph{Lefschetz fibrations and closed exotic $4$-manifolds}, minicourse at the Summer School on Modern Tools in Low-Dimensional Topology, ICTP, Trieste, June 2--6, 2025. School programme: \url{https://indico.ictp.it/event/10839/other-view?view=ictptimetable}.
\bibitem{Simons} Z.~Szab\'o, \emph{Exotic structures on smooth four-manifolds}, lecture at the Simons Collaboration on New Structures in Low-Dimensional Topology Annual Meeting, New York, March 27, 2026. Meeting programme and slides: \url{https://www.simonsfoundation.org/event/simons-collaboration-on-new-structures-in-low-dimensional-topology-annual-meeting-2026/}.
\bibitem{Xiao} G.~Xiao, \emph{Surfaces fibr\'ees en courbes de genre deux}, Lecture Notes in Math. 1137, Springer-Verlag, Berlin, 1985.
\end{thebibliography}
\end{document}